\documentclass[12pt, a4paper]{amsart}
\usepackage{amsthm, amssymb, color}
\usepackage{amsfonts, amscd}
\usepackage{epsfig}

\usepackage{xurl} 
\usepackage[colorlinks=true,  linkcolor=blue,  citecolor=green]{hyperref}

\usepackage{comment}

\usepackage{tikz}
\usepackage{tikz-cd}
\usepackage{tikz-3dplot}
\usetikzlibrary{shapes.geometric, calc}

\usetikzlibrary{math}

\theoremstyle{plain} \numberwithin{equation}{section}
\newtheorem{theorem}{Theorem}[section]
\newtheorem{corollary}[theorem]{Corollary}
\newtheorem{proposition}[theorem]{Proposition}
\newtheorem{lemma}[theorem]{Lemma}

\theoremstyle{definition}
\newtheorem{definition}[theorem]{Definition}
\newtheorem{example}[theorem]{Example}
\newtheorem{remark}[theorem]{Remark}
\newtheorem{construction}[theorem]{Construction}
\newtheorem{condition}[theorem]{Condition}

\def\Z{\mathbb Z}
\def\C{\mathbb C}
\def\R{\mathbb R}

\DeclareMathOperator{\rank}{rank}

\DeclareMathOperator{\Hom}{Hom}

\newcommand{\id}{\operatorname{id}}

\newcommand{\Th}{Th(\mathbb{C}P^1 \times \mathbb{C})}
\newcommand{\HT}{H^{\ast}_{T}}

\begin{document}

\title[Equivariant intersection cohomology of $T$-pseudomanifolds]
{Equivariant intersection cohomology of some $4$-dimensional locally standard
$T$-pseudomanifolds}

\author[Y. Koike]{Yuya Koike}
\address{Department of Pure and Applied Mathematics, Graduate School of Information Science and Technology, The University of Osaka, Suita, Osaka 565-0871, Japan.}
\email{u175137d@ecs.osaka-u.ac.jp}

\author[S. Kuroki]{Shintar\^o Kuroki}
\address{Department of Applied Mathematics, Faculty of Science, Okayama University of Science, 1-1 Ridai-Cho Kita-Ku Okayama-shi, Okayama 700-0005, Japan.}
\email{kuroki@ous.ac.jp}

\date{\today}

\begin{abstract}
Locally standard $T$-pseudomanifolds are singular analogues of quasitoric manifolds.
In this paper, we study a certain class of four-dimensional locally standard $T^2$-pseudomanifolds over convex polygons.
We establish an Orlik--Raymond type classification and study their ordinary and equivariant cohomology, as well as their intersection and equivariant intersection cohomology.
We also introduce ring structures on middle-perversity intersection cohomology and equivariant intersection cohomology via change of perversity.
We show that these rings can be described in terms of the cohomology and equivariant cohomology of associated nonsingular spaces, which are quasitoric manifolds or $S^4$.
Furthermore, we give combinatorial descriptions of the equivariant cohomology and equivariant intersection cohomology in terms of labelled graphs.
\end{abstract}

\maketitle 
\pagestyle{headings}
\tableofcontents

\section{Introduction}\label{sec:introduction}

A fundamental problem in toric topology is to understand the topology of spaces with torus actions in terms of the combinatorics of their orbit spaces together with isotropy data.
Quasitoric manifolds, introduced by Davis and Januszkiewicz, are locally standard $T$-manifolds whose orbit spaces are simple convex polytopes.
Their equivariant homeomorphism types are determined by characteristic pairs consisting of the orbit polytopes and characteristic functions.
The cohomology rings and equivariant cohomology rings of quasitoric manifolds can also be described in terms of the corresponding characteristic data \cite{DJ91,BP15}.

In our previous work \cite{KK25}, we introduced locally standard $T$-pseudomanifolds as singular analogues of quasitoric manifolds.
In this setting, the orbit space is allowed to be an arbitrary convex polytope, rather than a simple one.
The resulting class includes projective toric varieties viewed $T$-equivariantly.
We defined a characteristic functor for locally standard $T$-pseudomanifolds and, under suitable assumptions, established a classification in terms of characteristic data.
In \cite{Koi26}, the first author obtained a classification of all locally standard $T$-pseudomanifolds up to dimension three.
The present paper studies the topology and equivariant intersection cohomology of a class of four-dimensional locally standard $T$-pseudomanifolds over convex polygons.

\subsection{Main results}\label{subsec:main-results}
We now state the main results of the paper.
Throughout this paper, all cohomology groups are taken with coefficients in $\mathbb Z$.
Let $T=T^2$, and let $X$ be an oriented $4$-dimensional locally standard $T$-pseudomanifold whose orbit space $P=X/T$ is a convex polygon.
Let $F_1,\ldots,F_g$ be the edges of $P$, listed in cyclic order, and choose a primitive characteristic vector $\lambda_i\in\mathbb Z^2$ for each $F_i$.
Indices are understood modulo $g$.
We assume that
\[
\left|\det(\lambda_i,\lambda_{i+1})\right|\in\{0,1\}
\]
for every $i$, and that at least one determinant is nonzero.

A basic singular space arising in our setting is the Thom space of the trivial complex line bundle over $\mathbb CP^1$, $Th(\mathbb CP^1\times\mathbb C)$.
Its singular point has link $\mathbb CP^1\times S^1$
and is not an orbifold singularity.
In this paper, we allow such singularities
while excluding orbifold singularities.
Integral cohomology rings of $4$-dimensional toric orbifolds have been studied, for example, by Fu, So, and Song \cite{FSS23}.

Our first main result is an Orlik--Raymond type classification
that splits off Thom blocks corresponding to the singular points
(Theorem~\ref{thm:Orlik--Raymond-type-classification}).
It extends the connected-sum decomposition of
$4$-dimensional quasitoric manifolds \cite{OR70,BP02}
by allowing Thom blocks as additional building blocks.
Here a Thom block is weakly equivariantly homeomorphic to
$Th(\mathbb CP^1\times\mathbb C)$.

\begin{theorem}[Orlik--Raymond type classification]
\label{thm:intro-Orlik--Raymond-classification}
The space $X$ is equivariantly homeomorphic to an equivariant connected
sum of copies of $\mathbb CP^2$, Hirzebruch surfaces, and Thom blocks,
where the connected sums are taken at compatible regular fixed points.
\end{theorem}

Upon forgetting the torus action, this classification gives
\[
X\cong M\# B_1\#\cdots\# B_s,
\]
where $M$ is a quasitoric manifold or $S^4$,
each $B_j$ is a Thom block (possibly with reversed
orientation), and $s$ denotes the number of
singular points.

We next turn to equivariant intersection cohomology.

For singular spaces, ordinary cohomology does not in general satisfy Poincar\'e duality.
Intersection cohomology was introduced by Goresky and MacPherson to extend Poincar\'e duality to singular spaces \cite{GM80,GM83}.
We mainly follow the formulation in \cite{Fri20}.
We define equivariant intersection cohomology using finite-dimensional approximations to the Borel construction, following Kirwan \cite{Kir88}.
For the spaces considered in this paper, we define
a ring structure on middle-perversity
equivariant intersection cohomology.

\begin{theorem}[equivariant intersection cohomology]
\label{thm:intro-equivariant-intersection-cohomology}
There is an isomorphism of graded $H^\ast(BT)$-algebras
between the middle-perversity equivariant intersection
cohomology of $X$ and the equivariant cohomology of $M$:
\[
I_{\bar m}H_T^\ast(X)
\cong
H_T^\ast(M).
\]
\end{theorem}

The relation between intersection cohomology of toric
varieties and combinatorics has been studied by,
for example, Stanley \cite{Sta87} and Fieseler \cite{Fie91}.
For equivariant intersection cohomology, Barthel, Brasselet,
Fieseler, and Kaup gave a description in terms of sheaves
on fans \cite{BBFK99}.
Braden and MacPherson used moment graphs for certain
algebraic varieties \cite{BM01}.

By GKM theory \cite{GKM98,GZ01} and the theory
of torus graphs \cite{MMP07}, $H_T^\ast(M)$ can be
described by graph equivariant cohomology.
Thus, Theorem~\ref{thm:intro-equivariant-intersection-cohomology}
gives a combinatorial description of middle-perversity
equivariant intersection cohomology.

\subsection{Structure of the paper}
\label{subsec:structure-of-paper}

Sections~\ref{sec:locally-standard-T-pseudomanifold}
and~\ref{sec:canonical-models} recall locally standard
$T$-pseudomanifolds, characteristic functors,
and canonical models.
Section~\ref{sec:setting} describes the setting of the paper.

Part~\ref{part-1} is devoted to connected sums.
The Orlik--Raymond type classification is given in
Section~\ref{sec:decomposition-and-classification}.
The corresponding classification for $4$-dimensional
quasitoric manifolds is discussed, for example, in
\cite{BP02}.
A related Orlik--Raymond type classification for certain
$6$-dimensional torus manifolds was obtained by the
second author \cite{Kur16}.

Part~\ref{part-2} contains the cohomological computations.
Sections~\ref{sec:some-remarks}--\ref{sec:intersection-cohomology-rings}
describe the ordinary cohomology ring and intersection
cohomology of $X$ using the connected-sum decomposition.
Section~\ref{sec:equivariant-cohomology-labelled-graph}
gives a graph description of $H_T^\ast(X)$.
Section~\ref{sec:equivariant-intersection-cohomology}
studies the $H^\ast(BT)$-module and quotient ring structures
of equivariant intersection cohomology.
Finally, Section~\ref{sec:combinatorial-description-equivariant-intersection-cohomology}
proves Theorem~\ref{thm:intro-equivariant-intersection-cohomology}
using these results.
The Appendix, Section~\ref{sec:appendix-cohomology},
collects connected-sum formulas for ordinary and
intersection cohomology used in these computations.

\section{Locally standard \texorpdfstring{$T$}{T}-pseudomanifolds over convex polytopes}
\label{sec:locally-standard-T-pseudomanifold}

We recall locally standard $T$-pseudomanifolds introduced in
\cite{KK25}, restricting to the case where the orbit space is a convex
polytope.

Let $T=T^m$, and let $P$ be a $d$-dimensional convex polytope with
$d\leq m$.
For a compact space $L$, we write
\[
\mathring c(L):=L\times[0,1)/(L\times\{0\})
\]
for the open cone on $L$, and regard $\mathring c(\emptyset)$ as a
point.

We recall the definition of a locally standard $T$-pseudomanifold in the present setting.

\begin{definition}[locally standard $T$-pseudomanifold]
\label{def:locally-standard-T-pseudomanifold}
An $(m+d)$-dimensional \emph{locally standard $T$-pseudomanifold over
$P$} is a second-countable compact Hausdorff space $X$ equipped with a
continuous $T$-action such that $X/T=P$.
It is defined inductively on $d$ as follows.

If $d=0$, then $X=T$ with the multiplication action.
Suppose that $d>0$, and let $\pi:X\to P$ be the orbit projection.
For every $i$-dimensional face $F$ of $P$ and every
$x\in\pi^{-1}(\mathring F)$, we require the following:
\begin{enumerate}
    \item
    the isotropy subgroup $T_x$ is a $(d-i)$-dimensional subtorus of $T$;

\item
there exist a $T$-invariant open neighborhood $U_x$ of $x$ and a
space $L_x$.
If $i<d$, then $L_x$ is a $(2d-2i-1)$-dimensional locally standard
$T_x$-pseudomanifold over a $(d-i-1)$-dimensional convex polytope.
If $i=d$, we set $L_x=\emptyset$.
There exists a weakly $T$-equivariant homeomorphism
\[
U_x\cong
(\C^\ast)^i\times T^{m-d}\times\mathring c(L_x)
\]
via the isomorphism
\[
T\cong T/T_x\times T_x
\cong T^{i+m-d}\times T_x,
\]
where $T^{m-d+i}$ acts by multiplication on
$(\C^\ast)^i\times T^{m-d}$, and $T_x$ acts on
$\mathring c(L_x)$ through its action on $L_x$.
\end{enumerate}
\end{definition}

\begin{remark}
Quasitoric manifolds \cite{DJ91,BP15} are special cases of locally standard $T$-pseudomanifolds.
In this case, $m=d$.
For $i<d$, one may take
$L_x=S^{2d-2i-1}\subset\mathbb C^{d-i}$
with the coordinatewise torus action.
Thus $\mathring{c}(L_x)\cong\mathbb C^{d-i}$.
\end{remark}

\begin{remark}
The definition above is the specialization of
\cite[Definition~3.2]{KK25} to the case $l=0$ and $n=d$, with orbit
space a convex polytope.
In this setting, the filtration by orbit dimension is determined by
the faces of $P$.
Indeed, if $P^{(i)}$ denotes the $i$-skeleton of $P$, then the
corresponding term of the filtration is $\pi^{-1}(P^{(i)})$.
For every $i$-dimensional face $F$ of $P$, the preimage of the relative interior of $F$,
$\pi^{-1}(\mathring F)$,
is a topological manifold of dimension $m-d+2i$, and $\pi^{-1}(\mathring P)$ is the dense free stratum.
Thus the conditions on the orbit filtration appearing in
\cite{KK25} are automatic in the present setting.
\end{remark}

\begin{remark}
Since a convex polytope is a normal stratified pseudomanifold,
every locally standard $T$-pseudomanifold considered in this paper is
normal (see \cite{Koi26}).
Here, normal means that every link is connected.
\end{remark}

Let $\mathcal F(P)$ denote the face poset of $P$, regarded as a
category, and let $\mathcal T$ denote the category of subtori of $T$
with inclusions as morphisms.

We next describe the characteristic data associated with the face structure of $P$.

\begin{definition}[characteristic pair]\label{def:characteristic-pair}
A \emph{characteristic functor} on $P$ is a functor
\[
\lambda:\mathcal F(P)^{\mathrm{op}}\to\mathcal T
\]
such that
\[
\dim\lambda(F)=\operatorname{codim}F
\]
for every face $F$ of $P$.
The pair $(P,\lambda)$ is called a \emph{characteristic pair}.
\end{definition}

Let $X$ be a locally standard $T$-pseudomanifold over $P$.
For each face $F$ of $P$, the isotropy subgroup is constant on
$\pi^{-1}(\mathring F)$ (see \cite{KK25}).
Hence, for $x\in\pi^{-1}(\mathring F)$, the assignment
\[
\lambda(F):=T_x
\]
defines a characteristic functor on $P$.
We call $(P,\lambda)$ the \emph{characteristic pair of $X$}.

We also recall the notions of isomorphism and weak isomorphism for characteristic pairs.

\begin{definition}[isomorphisms of characteristic pairs]
Let $(P,\lambda)$ and $(P',\lambda')$ be characteristic pairs, and let
$f:P\to P'$ be a face-preserving homeomorphism.
We call $f$ a \emph{weak isomorphism} if there exists an automorphism
$\psi:T\to T$ such that $\lambda'(f(F))=\psi(\lambda(F))$ for every face $F$ of $P$.
If $\psi$ can be taken to be the identity, then $f$ is called an
\emph{isomorphism}.
\end{definition}

\begin{remark}
If $f=\id_P$, then $\lambda'$ is obtained from $\lambda$ by applying a single automorphism of $T$ to all characteristic subgroups.
\end{remark}

The classification theorem of \cite{KK25} specializes in the present
setting as follows.

\begin{theorem}[{\cite[Theorem~11.1]{KK25}}]
There is a one-to-one correspondence
\[
\frac{
\left\{
\begin{array}{c}
\text{locally standard $T$-pseudomanifolds}\\
\text{over convex polytopes}
\end{array}
\right\}
}{
\text{$T$-equivariant homeomorphism}
}
\quad\longleftrightarrow\quad
\frac{
\left\{
\begin{array}{c}
\text{characteristic pairs}\\
(P,\lambda)
\end{array}
\right\}
}{
\text{isomorphism}
}.
\]
Similarly, weakly $T$-equivariant homeomorphism classes correspond to
weak isomorphism classes of characteristic pairs.
\end{theorem}

\section{Canonical models}
\label{sec:canonical-models}

In this section, we introduce the canonical model associated with a characteristic pair and recall its relation to locally standard $T$-pseudomanifolds.
We first give the construction of the canonical model.

\begin{definition}[canonical model]\label{def:canonical-model}
Let $(P,\lambda)$ be a characteristic pair.
For $p\in P$, let $F(p)$ denote the unique face of $P$ whose relative
interior $\mathring{F(p)}$ contains $p$.
Define an equivalence relation $\sim_\lambda$ on $P\times T$ by
\[
(p,t)\sim_\lambda(p,t')
\quad\Longleftrightarrow\quad
t^{-1}t'\in\lambda(F(p)).
\]
The \emph{canonical model associated with $(P,\lambda)$} is
\[
X(P,\lambda):=(P\times T)/{\sim_\lambda}.
\]
The multiplication action of $T$ on the second factor induces a
$T$-action on $X(P,\lambda)$.
\end{definition}

The classification theorem identifies every locally standard $T$-pseudomanifold over a convex polytope with its canonical model.

\begin{theorem}\label{thm:canonical-model}
Let $X$ be a locally standard $T$-pseudomanifold over $P$, and let
$(P,\lambda)$ be its characteristic pair.
Then there exists a $T$-equivariant homeomorphism
\[
X\cong X(P,\lambda).
\]
\end{theorem}

\begin{remark}\label{rem:orientations}
Suppose that $m=d$, and fix orientations of $T$ and $P$.
We orient $X(P,\lambda)$ by the product orientation on
$T\times\mathring P$ via $(t,p)\mapsto[p,t]$.
We give $X$ the corresponding orientation under the
identification over $P$ in Theorem~\ref{thm:canonical-model}.
See \cite[Section~4]{BPR07} for the case of quasitoric manifolds.
\end{remark}

\section{The four-dimensional setting}\label{sec:setting}

The main results of this paper concern a particular class of $4$-dimensional locally standard $T^2$-pseudomanifolds.
In this section, we fix the assumptions and notation that will be used throughout the paper, unless otherwise stated.
In the section on connected sums
(Section~\ref{sec:connected-sum}),
some constructions will temporarily be formulated in arbitrary dimension.

Let $T=T^2$ and let $X$ be a $4$-dimensional locally standard
$T$-pseudomanifold whose orbit space $P=X/T$ is a convex polygon.
Namely, $m=d=2$ in Definition~\ref{def:locally-standard-T-pseudomanifold}.
We identify $P$ with a convex polygon in $\mathbb R^2$
and give it the standard orientation.
We orient $X$ as in Remark~\ref{rem:orientations}.
Let $\pi:X \to P$ denote the orbit projection.
Let $F_1,\ldots,F_g$ be the edges of $P$, listed in counterclockwise cyclic order with respect to this orientation, and set
\[
v_i:=F_i\cap F_{i+1},
\qquad
1\leq i\leq g,
\]
where indices are understood modulo $g$.
Let $\lambda:\mathcal F(P)^{\mathrm{op}}\to\mathcal T$
be the characteristic functor of $X$
(see Definition~\ref{def:characteristic-pair}).
For each edge $F_i$, set
\[
K_i:=\lambda(F_i)\subset T.
\]
Then $K_i$ is a circle subgroup of $T$, while $\lambda(v_i)=T^2$ for every vertex $v_i$.
Let $N:=\Hom(S^1,T)\cong\mathbb Z^2$ be the lattice of one-parameter subgroups of $T$.
For each $i$, choose a primitive vector $\lambda_i\in N$, where $\exp(\mathbb R\lambda_i)=K_i$.
Such a vector is uniquely determined up to sign.

We use the following terminology for the chosen primitive vectors.

\begin{definition}[characteristic vector]\label{def:characteristic-vector}
A primitive vector $\lambda_i\in N$ representing the circle subgroup $K_i=\lambda(F_i)$ is called a \emph{characteristic vector} of the edge $F_i$.
\end{definition}

The characteristic vectors determine an integer associated with each vertex.

\begin{definition}[vertex index]\label{def:vertex-index}
For a vertex $v_i=F_i\cap F_{i+1}$,
its \emph{index} is defined by
\[
\operatorname{ind}(v_i)
:=
\left|
\det(\lambda_i,\lambda_{i+1})
\right|.
\]
\end{definition}

We distinguish the following two types of vertices according to their indices.

\begin{definition}[regular vertices and $0$-vertices]\label{def:name-of-special-vertex}
A vertex $v$ is called a \emph{regular vertex} if $\operatorname{ind}(v)=1$.
A vertex $v$ is called a \emph{$0$-vertex} if $\operatorname{ind}(v)=0$.
\end{definition}

We now state the condition on vertex indices used throughout the paper.

\begin{condition}[0-1 condition]\label{cond:01-condition}
Every vertex of $P$ is either a regular vertex or a $0$-vertex; equivalently,
\[
\operatorname{ind}(v_i)\in\{0,1\}
\qquad
\text{for every }i.
\]
Moreover, we assume that $P$ has at least one regular vertex.
\end{condition}

Let $s$ and $r$ denote the numbers of $0$-vertices and regular vertices of $P$, respectively.
Then $g=r+s$.

\begin{remark}\label{rem:at-least-two-regular-vertices}
Under Condition~\ref{cond:01-condition}, $P$ has at least two regular vertices.
Indeed, if there were exactly one regular vertex, then all the other vertices would be $0$-vertices. 
Since the characteristic vectors on the two edges meeting at a $0$-vertex agree up to sign, going around the boundary of $P$ shows that all characteristic vectors agree up to sign. 
This would force the remaining vertex to be a $0$-vertex as well, a contradiction.
\end{remark}

Thus, we have $r\geq2$.

\part{Connected sums and the Orlik--Raymond type classification}\label{part-1}

In this part, we establish an Orlik--Raymond type classification of the spaces considered in this paper. For this purpose, we introduce Thom blocks and develop equivariant connected sums of characteristic pairs and their canonical models. The resulting classification, Theorem~\ref{thm:Orlik--Raymond-type-classification}, gives the connected-sum decomposition in Corollary~\ref{cor:decomposition-M-Thom}.
We combine this decomposition with the connected-sum formulas in Theorems~\ref{thm:connected-sum-for-cohomology-ring} and~\ref{thm:connected-sum-for-intersection-cohomology} in the Appendix.
This reduces the computations in Part~\ref{part-2} to those for a quasitoric manifold or $S^4$ and for the Thom blocks.

\section{Thom triangle and Thom block}
\label{sec:Thom-triangle-block}

In this section, we study characteristic pairs over a triangle satisfying
Condition~\ref{cond:01-condition}.
We show that every such characteristic pair whose canonical model is not a quasitoric manifold is a \emph{Thom triangle}.
We then define the canonical model associated with a Thom triangle to be
a \emph{Thom block} and show that every Thom block is weakly
$T^2$-equivariantly homeomorphic to the Thom space
$Th(\mathbb CP^1\times\mathbb C)$.
Thom blocks will serve as basic building blocks in the
Orlik--Raymond type classification developed in the following sections.

We first classify the characteristic pairs over a triangle satisfying Condition~\ref{cond:01-condition}.

\begin{proposition}[characteristic pairs over a triangle]
\label{prop:characteristic-functor-on-triangle}
Let $(P,\lambda)$ be a characteristic pair satisfying Condition~\ref{cond:01-condition}, where $P$ is a triangle.
Then exactly one of the following holds.
\begin{enumerate}
\item Every vertex $v$ of $P$ satisfies $\operatorname{ind}(v)=1$.
In this case, $X(P,\lambda)$ is a quasitoric manifold over a triangle;
in particular, it is weakly equivariantly homeomorphic to $\mathbb CP^2$.

\item Exactly one vertex $v$ of $P$ satisfies $\operatorname{ind}(v)=0$,
and the other two vertices have index $1$.
Moreover, after applying an automorphism of $T^2$, the characteristic subgroups on the three edges are
\[
T\langle(1,0)\rangle,\qquad
T\langle(0,1)\rangle,\qquad
T\langle(1,0)\rangle.
\]
\end{enumerate}
\end{proposition}

\begin{proof}
By Remark~\ref{rem:at-least-two-regular-vertices}, at least two vertices of $P$ have index $1$.
Since $P$ has three vertices, either all three vertices have index $1$, or exactly one has index $0$.
In the former case, it is well known that $X(P,\lambda)$ is weakly $T^2$-equivariantly homeomorphic to $\mathbb CP^2$.
It remains to consider the latter case.
Choose primitive characteristic vectors $\lambda_1,\lambda_2,\lambda_3\in\mathbb Z^2$ as in Figure~\ref{fig:char-functor-on-triangle}.
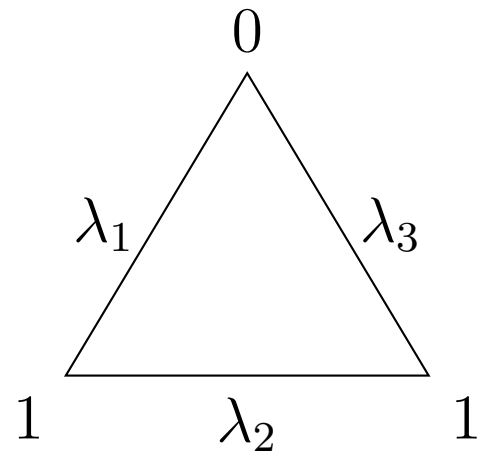
\begin{figure}[htbp]
\centering
\begin{tikzpicture}[scale=0.5]
\coordinate (A) at (0,2.8);
\coordinate (B) at (-2.4,-1.2);
\coordinate (C) at (2.4,-1.2);

\draw (A)--(B)--(C)--cycle;

\node[left] at ($(A)!0.5!(B)$) {$\lambda_1$};
\node[below] at ($(B)!0.5!(C)$) {$\lambda_2$};
\node[right] at ($(C)!0.5!(A)$) {$\lambda_3$};

\node[above] at (A) {$0$};
\node[below left] at (B) {$1$};
\node[below right] at (C) {$1$};
\end{tikzpicture}
\caption{Characteristic vectors in the case with one vertex of index $0$. The numbers at the vertices indicate the vertex indices.}
\label{fig:char-functor-on-triangle}
\end{figure}
Since $\lambda_1$ and $\lambda_2$ form a basis of $\mathbb Z^2$, after applying an automorphism of $T^2$, we may assume that $\lambda_1=(1,0)$ and $\lambda_2=(0,1)$.
Then, because the remaining vertex is a $0$-vertex, $\lambda_3$ and $\lambda_1$ are dependent.
Since $\lambda_3$ is primitive, $\lambda_3=(1,0)$ or $(-1,0)$.
Therefore, the corresponding characteristic subgroups are $T\langle(1,0)\rangle, T\langle(0,1)\rangle, T\langle(1,0)\rangle$.
\end{proof}

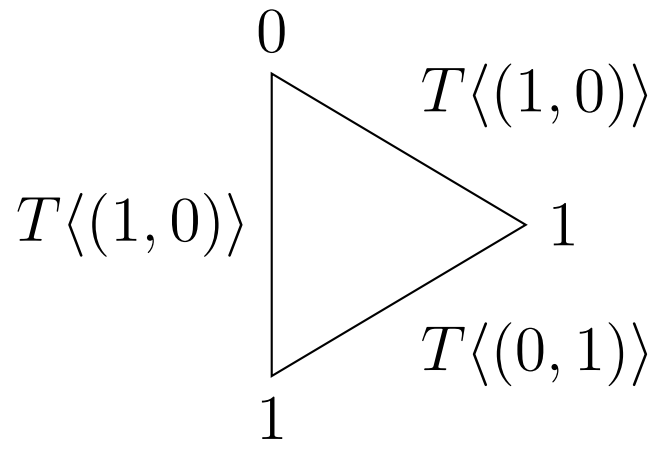
\begin{figure}[htbp]
\centering
\begin{tikzpicture}[scale=0.08]
\coordinate (A) at (0,30);
\coordinate (B) at (0,5);
\coordinate (C) at (21,17.5);

\draw (A)--(B)--(C)--cycle;

\node[left] at ($(A)!0.5!(B)$) {$T\langle(1,0)\rangle$};
\node[below right] at ($(B)!0.5!(C)$) {$T\langle(0,1)\rangle$};
\node[above right] at ($(C)!0.5!(A)$) {$T\langle(1,0)\rangle$};

\node[above] at (A) {$0$};
\node[below] at (B) {$1$};
\node[right] at (C) {$1$};
\end{tikzpicture}
\caption{A characteristic pair representing a Thom triangle. The numbers at the vertices indicate the vertex indices.}
\label{fig:Thom-triangle}
\end{figure}

The second case of Proposition~\ref{prop:characteristic-functor-on-triangle} motivates the following definition.

\begin{definition}[Thom triangle]\label{def:Thom-triangle}
A characteristic pair over a triangle is called a \emph{Thom triangle}
if it is weakly isomorphic to the characteristic pair shown in
Figure~\ref{fig:Thom-triangle}.
We denote a Thom triangle by
\[
(\Delta_{Th},\lambda_{Th}).
\]
\end{definition}

By Proposition~\ref{prop:characteristic-functor-on-triangle}, every characteristic pair over a triangle satisfying Condition~\ref{cond:01-condition} and having a $0$-vertex is a Thom triangle.

We next identify the canonical model associated with a Thom triangle.

\begin{proposition}\label{prop:canonical-model-of-Thom-triangle}
The canonical model associated with a Thom triangle is weakly $T^2$-equivariantly homeomorphic to the Thom space
$Th(\mathbb CP^1\times\mathbb C)$.
\end{proposition}

\begin{proof}
By \cite[Example~7.1]{KK25}, the Thom space
$Th(\mathbb CP^1\times\mathbb C)$ has the characteristic pair shown in Figure~\ref{fig:Thom-triangle} as its characteristic pair.
Hence the assertion follows from the classification theorem
\cite[Theorem~11.1]{KK25}.
\end{proof}

In view of Proposition~\ref{prop:canonical-model-of-Thom-triangle}, we introduce the following terminology.

\begin{definition}[Thom block]\label{def:Thom-block}
The canonical model associated with a Thom triangle is called a
\emph{Thom block}.
\end{definition}

\section{Connected sums}\label{sec:connected-sum}

In this section, we introduce connected sums of characteristic pairs and show that the corresponding canonical model is equivariantly homeomorphic to the equivariant connected sum of the canonical models (Theorem~\ref{thm:connected-sum}).
Although the main results of this paper concern the four-dimensional case, we work temporarily in arbitrary dimension.
Thus, throughout this section, $P$ and $Q$ denote $d$-dimensional convex polytopes with $d\geq 2$, and $T=T^d$.

\subsection{Connected sums of convex polytopes}
\label{subsec:connected-sum-of-polytope}

We introduce a connected sum for convex polytopes at simple vertices, modeled on the connected-sum for simple polytopes in \cite[Construction~1.1.13]{BP15} and \cite[Section~6]{BR01}.

We first introduce the notion of a simple vertex.

\begin{definition}[Simple vertex]\label{def:simple-vertex}
A vertex $v$ of a $d$-dimensional convex polytope is called \emph{simple} if exactly $d$ facets contain $v$.
\end{definition}

We now define the connected sum of two convex polytopes at simple vertices.

\begin{construction}[connected sum of convex polytopes]
\label{cons:connected-sum-of-polytope}
Let $P$ and $Q$ be $d$-dimensional convex polytopes with $d\geq2$,
and let $v\in P$ and $w\in Q$ be simple vertices.
Fix a pairing between the facets containing $v$ and those containing $w$.
We describe the construction informally.
Truncate $P$ at $v$ and $Q$ at $w$.
Since $v$ and $w$ are simple, the truncations produce new facets
$\Delta(v)$ and $\Delta(w)$, each of which is a $(d-1)$-simplex.
We then glue the two truncated polytopes by identifying
$\Delta(v)$ with $\Delta(w)$ in such a way that each pair of
corresponding facets incident to $v$ and $w$ combines to form
a single facet.
The gluing is carried out by suitable projective transformations.
For precise constructions, see
\cite[Construction~1.1.13]{BP15}
and \cite[Section~6]{BR01},
where the resulting space is defined as a
$d$-dimensional convex polytope.
We note from these constructions that the simplicity assumption
is needed only at $v$ and $w$.

We define this polytope to be the connected sum of $P$ and $Q$
at $v$ and $w$, and denote it by $P\#_{v,w}Q$.
When the vertices and the pairing of the incident facets are clear
from the context, we simply write $P\#Q$.
\end{construction}

\begin{remark}
The combinatorial type of $P\#_{v,w}Q$ may depend on the chosen pairing between the facets containing $v$ and those containing $w$.
\end{remark}

\subsection{Connected sums of characteristic pairs}
\label{subsec:connected-sum-of-char-pair}

We now extend the connected-sum construction to characteristic pairs.
We first introduce the sign of a regular vertex, which will be used in the connected-sum construction.
Our convention is modeled on the vertex signs used for omnioriented quasitoric manifolds
(see \cite[Section~5]{BPR07}).
Let $(P,\lambda)$ be a characteristic pair, where $P \subset \R^d$ is an oriented $d$-dimensional convex polytope.
For each facet $F$ of $P$, the characteristic functor $\lambda$ assigns a circle subgroup $\lambda(F)\subset T$.
This subgroup is generated by a primitive vector in $\mathbb Z^d$, uniquely determined up to sign.
For each facet $F$, choose one of the two primitive vectors and denote it by $\lambda_F$.

\begin{remark}
In the case of quasitoric manifolds, a choice of an orientation of the orbit polytope together with a choice of sign for each characteristic vector is equivalent to a choice of an omniorientation.
Thus, the choices of characteristic vectors above may be regarded as the analogue of the sign data determined by an omniorientation.
\end{remark}

Let $v$ be a simple vertex of $P$, and let $F_{v,1},\ldots,F_{v,d}$ be the facets containing $v$.
For each $j$, let $E_{v,j}$ be the edge at $v$ opposite to $F_{v,j}$, and choose a vector $e_{v,j}$ pointing from $v$ along $E_{v,j}$.
We order $F_{v,1},\ldots,F_{v,d}$ so that
$(e_{v,1},\ldots,e_{v,d})$ is positively oriented with respect to the orientation of $P$.

The chosen characteristic vectors at a simple vertex form the following matrix.

\begin{definition}[vertex matrix]\label{def:vertex-matrix}
The \emph{vertex matrix} at $v$ is
\[
\Lambda_v:=
(\lambda_{F_{v,1}},\ldots,\lambda_{F_{v,d}}).
\]
\end{definition}

The determinant $\det\Lambda_v$ depends on the choices
of signs of the characteristic vectors
and on the orientation of $P$.
We use the determinant of the vertex matrix to define regular vertices, $0$-vertices, and vertex signs.

\begin{definition}[regular vertex, $0$-vertex, vertex sign]
\label{def:vertex-sign}
A simple vertex $v$ is called a \emph{regular vertex} if $|\det\Lambda_v|=1$, and a \emph{$0$-vertex} if $\det\Lambda_v=0$.
For a regular vertex $v$, its \emph{vertex sign} is defined by
\[
\sigma(v):=\det\Lambda_v\in\{1,-1\}.
\]
\end{definition}

\begin{remark}
  When $d=2$, the notions of regular vertex and $0$-vertex agree with those introduced in Section~\ref{sec:setting}.
\end{remark}

We now consider two characteristic pairs $(P,\lambda)$ and $(Q,\mu)$,
where $P$ and $Q$ are oriented $d$-dimensional convex polytopes.

\begin{condition}[Compatibility Condition]
\label{cond:compatibility}
Let $v\in P$ and $w\in Q$ be regular vertices.
We call $v$ and $w$ compatible if the facets containing them admit a pairing satisfying the following conditions:
\begin{enumerate}
\item
For every pair of corresponding facets $F$ and $G$,
we have $\lambda(F)=\mu(G)$.
\item
Choose primitive characteristic vectors so that $\lambda_F=\mu_G$ on corresponding facets.
We require that the vertex signs satisfy $\sigma(v)=-\sigma(w)$.
\end{enumerate}
\end{condition}

\begin{remark}
The Compatibility Condition is modeled on the sign condition for connected sums of omnioriented quasitoric manifolds.
In particular, \cite[Proposition~5.3]{BPR07} shows that the corresponding connected sum is compatible with the given orientations precisely when the two vertex signs are opposite.
\end{remark}

We can therefore glue the two characteristic functors
along the corresponding faces.

\begin{construction}[connected sum of characteristic pairs]
\label{cons:connected-sum-of-characteristic-pairs}
Let $(P,\lambda)$ and $(Q,\mu)$ be characteristic pairs, where $P$ and $Q$ are oriented, and let $v\in P$ and $w\in Q$ be compatible regular vertices.
Using the pairing in the Compatibility Condition, form the connected sum $P\#_{v,w}Q$ as in Construction~\ref{cons:connected-sum-of-polytope}.
We define a characteristic functor $\lambda \#\mu$ on $P\#_{v,w}Q$ as follows.

A face inherited from the $P$-part is assigned the same characteristic subgroup as the corresponding face of $P$, and a face inherited from the $Q$-part is assigned the subgroup given by $\mu$.
If a face $H$ of $P\#_{v,w}Q$ is obtained by joining a face $F$ of $P$ containing $v$ with the corresponding face $G$ of $Q$ containing $w$, we set
\[
(\lambda\#\mu)(H):=\lambda(F)=\mu(G).
\]
We call $(P\#_{v,w}Q,\lambda \#\mu)$ the \emph{connected sum of $(P,\lambda)$ and $(Q,\mu)$ at $v$ and $w$}.
When the vertices are clear from the context, we simply write $(P,\lambda) \#(Q,\mu)$.
\end{construction}

The equality $\lambda(F)=\mu(G)$ for corresponding faces shows that $\lambda\#\mu$ is well-defined and is a characteristic functor on $P\#_{v,w}Q$.

\begin{remark}
The Compatibility Condition ensures that the identification
of the truncation facets $\Delta(v)$ and $\Delta(w)$
reverses their induced boundary orientations.
Hence $P\#_{v,w}Q$ inherits an orientation agreeing with
the given orientations on both truncated polytopes.
\end{remark}

\begin{remark}\label{rem:equivariant-connected-sum}
Let $X=X(P,\lambda)$ and $Y=X(Q,\mu)$ be canonical models.
Let $x$ and $y$ be the fixed points over compatible regular
vertices.
The connected-sum construction of \cite[Definition~10.6]{BJ82}
can be performed $T$-equivariantly.
We denote the resulting space by $X\#_{x,y}Y$.
When $x$ and $y$ are clear from the context, we simply write $X\#Y$.
\end{remark}

\subsection{The equivariant connected-sum formula}
\label{subsec:equivariant-connected-sum-formula}

We now compare the equivariant connected sum above with the canonical
model associated with the connected sum of characteristic pairs.
The following theorem identifies this equivariant connected sum with the canonical model of the connected sum of characteristic pairs.

\begin{theorem}\label{thm:connected-sum}
Let $(P,\lambda)$ and $(Q,\mu)$ be characteristic pairs, and let
$v\in P$ and $w\in Q$ be regular vertices satisfying the
Compatibility Condition. Then there exists a $T$-equivariant
homeomorphism
\[
X(P,\lambda)\#X(Q,\mu)
\cong
X(P\# Q,\lambda\#\mu),
\]
where the left-hand side denotes the equivariant connected sum of the
canonical models $X(P,\lambda)$ and $X(Q,\mu)$ at the fixed points
corresponding to $v$ and $w$, and the right-hand side is the canonical
model associated with the connected sum $(P\# Q,\lambda\#\mu)$ of the
characteristic pairs.
\end{theorem}

\begin{proof}
Put $X:=X(P,\lambda)$, $Y:=X(Q,\mu)$, and
$Z:=X(P\#Q,\lambda\#\mu)$.
Since $v$ and $w$ are regular, the corresponding fixed points
have $T$-invariant neighborhoods $T$-equivariantly homeomorphic
to $\mathbb C^d$.
By the Compatibility Condition, we may use the same
linear $T$-action in both cases.
The orbit spaces of these neighborhoods correspond to
the regions cut off near $v$ and $w$ in $P$ and $Q$,
respectively.
The identification of the boundary spheres is induced
by the gluing of $\Delta(v)$ and $\Delta(w)$.
The embeddings of the truncated polytopes into $P\#Q$
induce continuous $T$-equivariant maps
from the preimages of these polytopes in $X$ and $Y$
to $Z$, given by
\[
[p,t]\longmapsto[p',t].
\]
Here $p$ is a point of the truncated polytope of $P$ or $Q$,
and $t\in T$.
The class $[p,t]$ is regarded as a point of $X$ or $Y$,
respectively, and $p'$ denotes the image of $p$ under
the corresponding embedding into $P\#Q$.
These maps are well-defined and injective because
the characteristic subgroups are preserved.
Moreover, the maps agree under the identification of
the boundary spheres, since the characteristic subgroups
agree on corresponding faces.
Therefore, they define a continuous $T$-equivariant map
$J:X\#Y\to Z$.
By construction, the preimages in $Z$ of the two
truncated polytopes cover $Z$ and
intersect precisely along the identified boundary.
Hence $J$ is bijective.
Since $X\#Y$ is compact and $Z$ is Hausdorff,
$J$ is a homeomorphism.
\end{proof}

\section{Decomposition and classification}
\label{sec:decomposition-and-classification}

In this section, we apply the connected-sum construction to successively split off Thom triangles from the characteristic pair.
We then derive the Orlik--Raymond type classification and the collapse description of $X$.

\subsection{Thom-triangle decomposition and classification}
\label{subsec:Thom-triangle-decomposition}

We first show that a Thom triangle can be split off whenever a $0$-vertex is present.
Let $(P,\lambda)$ be a characteristic pair satisfying Condition~\ref{cond:01-condition}, and suppose that $P$ is a $g$-gon with $g\geq4$.

\begin{lemma}[Thom-triangle decomposition]
\label{lem:Thom-triangle-decomposition}
Suppose that $s>0$.
Then there exists a characteristic pair $(P',\lambda')$ over a $(g-1)$-gon such that
\[
(P,\lambda)
\cong
(\Delta_{Th},\lambda_{Th})
\#(P',\lambda'),
\]
where $(\Delta_{Th},\lambda_{Th})$ is a Thom triangle.
Moreover, $(P',\lambda')$ satisfies Condition~\ref{cond:01-condition}, has $(s-1)$ $0$-vertices, and has $r$ regular vertices.
\end{lemma}

\begin{proof}
Since $P$ has both a $0$-vertex and a regular vertex, there exist adjacent vertices $v_0$ and $v_1$ such that $v_0$ is a $0$-vertex and $v_1$ is regular.
Let $F$ be the edge joining them, and let $F_0$ and $F_1$ be the other edges containing $v_0$ and $v_1$, respectively (see Figure~\ref{fig:char-pair-local}).

\begin{figure}[htbp]
\centering
\begin{tikzpicture}

\coordinate (A) at (-1,1.732);
\coordinate (B) at ( 1,1.732);
\coordinate (D) at ( 1,-1.732);
\coordinate (E) at (-1,-1.732);
\coordinate (V1) at (-2.0,-0.866);
\coordinate (V0) at (-2.0,0.866);

\draw (V0)--(A);
\draw (A)--(B);
\draw (V1)--(V0);

\draw[dashed] (B)--(1.7,1.0);
\draw[dashed] (1.7,-1.0)--(D);

\draw (D)--(E);
\draw (E)--(V1);

\node at (2.0,0) {$\vdots$};

\draw[-stealth] (-0.5,0)
 arc[start angle=-180,end angle=90,radius=0.5];

\node[left] at ($(V0)!0.5!(V1)$) {$\lambda_{F}$};
\node[right] at ($(V0)!0.5!(V1)$) {$F$};

\node[below left] at ($(E)!0.5!(V1)$) {$\lambda_{F_1}$};
\node[right] at ($(E)!0.5!(V1)$) {$F_1$};

\node[above left] at ($(V0)!0.5!(A)$) {$\lambda_{F_0}$};
\node[below right] at ($(V0)!0.5!(A)$) {$F_0$};

\node[left] at (V1) {$\pm 1$};
\node[above right] at (V1) {$v_1$};

\node[left] at (V0) {$0$};
\node[below right] at (V0) {$v_0$};

\end{tikzpicture}

\caption{A local configuration of the characteristic pair $(P,\lambda)$. The signs $\pm1$ at regular vertices indicate the corresponding vertex signs.}
\label{fig:char-pair-local}
\end{figure}
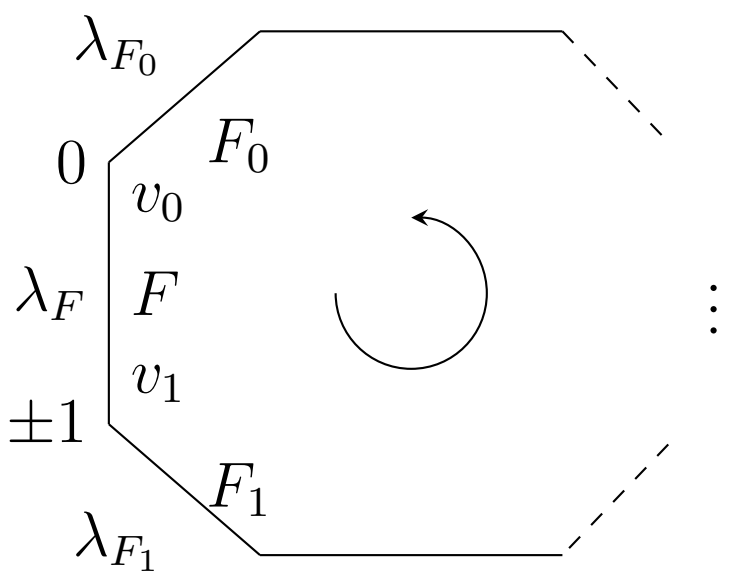
Since $v_0$ is a $0$-vertex, we have $\lambda_{F_0}=\pm\lambda_F$.
Hence we may remove $v_0$ and replace the two consecutive edges $F_0$ and $F$ by a single edge.
This produces a characteristic pair $(P',\lambda')$ over a $(g-1)$-gon, as shown in Figure~\ref{fig:Thom-triangle-decomposition}.

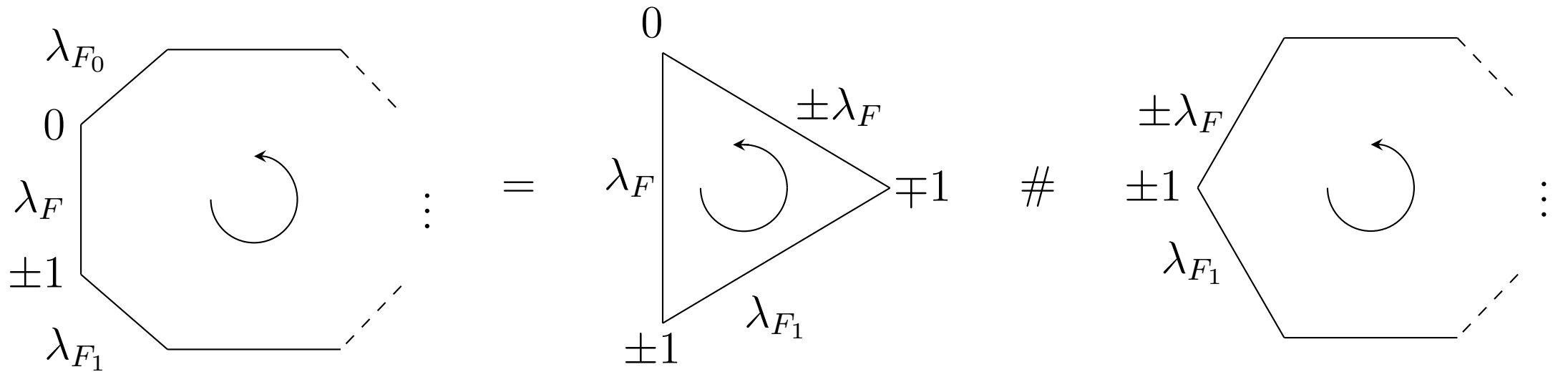
\begin{figure}[htbp]
\centering

\begin{tabular}{c c c c c}

\raisebox{-1.732cm}{
\begin{tikzpicture}[scale=0.8]
\coordinate (A) at (-1,1.732);
\coordinate (B) at ( 1,1.732);
\coordinate (D) at ( 1,-1.732);
\coordinate (E) at (-1,-1.732);
\coordinate (F) at (-2.0,-0.866);
\coordinate (G) at (-2.0,0.866);

\draw (G)--(A);
\draw (A)--(B);
\draw (F)--(G);

\draw[dashed] (B)--(1.7,1.0);
\draw[dashed] (1.7,-1.0)--(D);

\draw (D)--(E);
\draw (E)--(F);

\node at (2.0,0) {$\vdots$};

\draw[-stealth] (-0.5,0)
 arc[start angle=-180,end angle=90,radius=0.5];

\node[left] at ($(G)!0.5!(F)$) {$\lambda_{F}$};
\node[below left] at ($(E)!0.5!(F)$) {$\lambda_{F_1}$};
\node[above left] at ($(G)!0.5!(A)$) {$\lambda_{F_0}$};

\node[left] at (F) {$\pm 1$};
\node[left] at (G) {$0$};
\end{tikzpicture}
}

&
{$=$}
&

\raisebox{-1.7cm}{
\begin{tikzpicture}[scale=0.1]
\draw (0,5) -- (0,30);
\draw (0,5) -- (21,17.5);
\draw (0,30) -- (21,17.5);

\draw (3.5,17.5) arc[start angle=180,end angle=360,radius=4];
\draw (11.5,17.5) arc[start angle=0,end angle=90,radius=4];

\draw[-stealth] (8.0,21.5) -- (6.5,21.5);

\node at (-3,18.5) {$\lambda_{F}$};
\node at (16.3,25) {$\pm \lambda_{F}$};
\node[below] at (10.5,9) {$\lambda_{F_1}$};

\node[above] at (-1,30) {$0$};
\node at (-1,2.5) {$\pm 1$};
\node at (24,17.5) {$\mp 1$};
\end{tikzpicture}
}

&
{$\#$}
&

\raisebox{-1.3cm}{
\begin{tikzpicture}[scale=0.8]
\coordinate (A) at (-1,1.732);
\coordinate (B) at ( 1,1.732);
\coordinate (D) at ( 1,-1.732);
\coordinate (E) at (-1,-1.732);
\coordinate (F) at (-2,0);

\draw (F)--(A);
\draw (A)--(B);

\draw[dashed] (B)--(1.7,1.0);
\draw[dashed] (1.7,-1.0)--(D);

\draw (D)--(E);
\draw (E)--(F);

\node at (2.0,0) {$\vdots$};

\draw[-stealth] (-0.5,0)
 arc[start angle=-180,end angle=90,radius=0.5];

\node[left] at ($(E)!0.5!(F)$) {$\lambda_{F_1}$};
\node[left] at ($(F)!0.5!(A)$) {$\pm \lambda_{F}$};

\node[left] at (F) {$\pm1$};
\end{tikzpicture}
}

\end{tabular}

\caption{The connected-sum decomposition of the characteristic pair $(P,\lambda)$.}
\label{fig:Thom-triangle-decomposition}
\end{figure}

By construction, $(P,\lambda)$ is recovered from $(P',\lambda')$ by taking the connected sum with a Thom triangle at the regular vertex corresponding to $v_1$.
We orient the Thom triangle and $P'$ so that their
orientations agree with that of $P$ on the corresponding
truncated pieces.
The gluing vertices then satisfy the Compatibility Condition,
and we obtain
\[
(P,\lambda)\cong (\Delta_{Th},\lambda_{Th})\#(P',\lambda').
\]

The construction removes exactly one $0$-vertex and does not change any other vertex index.
Hence $(P',\lambda')$ satisfies Condition~\ref{cond:01-condition}, has $(s-1)$ $0$-vertices, and has $r$ regular vertices.
\end{proof}

We now iterate the Thom-block decomposition obtained above.
In the case of quasitoric manifolds of dimension $4$,
we use the formulation given in
\cite[Theorem~5.65 and Corollary~5.66]{BP02}.

\begin{theorem}[Orlik--Raymond type classification]
\label{thm:Orlik--Raymond-type-classification}
Let $X$ be a $4$-dimensional locally standard $T^2$-pseudomanifold satisfying Condition~\ref{cond:01-condition}.
Then $X$ is equivariantly homeomorphic to an equivariant connected sum of spaces of the following types:
\begin{itemize}
\item $\mathbb CP^2$;
\item Hirzebruch surfaces $H_k$;
\item Thom blocks.
\end{itemize}
\end{theorem}

\begin{proof}
Let $(P,\lambda)$ be the characteristic pair of $X$.
If $s=0$, then $X$ is a quasitoric manifold, and the assertion follows from \cite[Theorem~5.65]{BP02}.
Suppose that $s>0$.
Repeated application of
Lemma~\ref{lem:Thom-triangle-decomposition}
and Theorem~\ref{thm:connected-sum}
splits off Thom blocks while decreasing the number of
$0$-vertices and preserving the number of regular vertices.
If $r\geq3$, this process eventually leaves a quasitoric manifold, and the assertion follows from the classification result for $4$-dimensional quasitoric manifolds.
If $r=2$, it terminates with a triangle having one $0$-vertex and two regular vertices, which is a Thom triangle by Proposition~\ref{prop:characteristic-functor-on-triangle}.
The assertion follows.
\end{proof}

Forgetting the torus action gives the following nonequivariant decomposition.

\begin{corollary}
\label{cor:decomposition-M-Thom}
Let $s$ and $r$ be the numbers of $0$-vertices and regular vertices of $P$, respectively.
Then, forgetting the equivariant structure, we may write
\[
X\cong M\# B_1\#\cdots\# B_s,
\]
where each $B_j$ is either
$Th(\mathbb CP^1\times\mathbb C)$
or
$\overline{Th(\mathbb CP^1\times\mathbb C)}$,
and
\[
M=
\begin{cases}
\text{a quasitoric manifold over an $r$-gon}, & r\geq3,\\
S^4, & r=2.
\end{cases}
\]
Here $\overline{Y}$ denotes the space $Y$ with the reversed orientation.
The orientation of each Thom block is determined by the Compatibility Condition.
In the case $r=2$, the factor $S^4$ is introduced only as the neutral element for the ordinary connected sum.
\end{corollary}

\begin{proof}
This follows from the proof of Theorem~\ref{thm:Orlik--Raymond-type-classification}.
\end{proof}

\subsection{Resolution by a nonsingular space}
\label{subsec:resolution-by-quasitoric-manifold}

The decomposition above also gives another description of $X$ in terms of a quasitoric manifold (or $S^4$).
The following proposition expresses this description in terms of collapsing circle orbits.

\begin{proposition}[resolution by a quasitoric manifold or $S^4$]
\label{prop:resolution-by-quasitoric}
Let $X$ be a $4$-dimensional locally standard $T^2$-pseudomanifold satisfying Condition~\ref{cond:01-condition}, and let $r$ be the number of regular vertices of its orbit polygon.
Then $X$ is obtained by collapsing finitely many circle orbits of
\[
M=
\begin{cases}
\text{a quasitoric manifold over an $r$-gon}, & r\geq3,\\
S^4, & r=2,
\end{cases}
\]
to points.
\end{proposition}

\begin{proof}
For $r\geq3$, remove the $0$-vertices successively, replacing the two edges meeting at each $0$-vertex by a single edge carrying their common characteristic subgroup.
The resulting characteristic pair has only regular vertices and hence determines a quasitoric manifold over an $r$-gon.
For $r=2$, the corresponding nonsingular part is $S^4$ with its standard $T^2$-action.
Reversing the above operations amounts to collapsing the corresponding circle orbits to points.
\end{proof}

\part{Cohomology}\label{part-2}

In this part, we compute the ordinary, intersection, equivariant, and equivariant intersection cohomology of the spaces considered in this paper. Combining the connected-sum decomposition of Corollary~\ref{cor:decomposition-M-Thom} with the connected-sum formulas in Theorems~\ref{thm:connected-sum-for-cohomology-ring} and~\ref{thm:connected-sum-for-intersection-cohomology} of the Appendix, we determine the ordinary cohomology ring and the intersection cohomology groups in Theorems~\ref{thm:ordinary-cohomology-ring} and~\ref{thm:intersection-cohomology-group}, respectively. We then use change of perversity to define a ring structure on middle-perversity intersection cohomology and identify the resulting ring with $H^\ast(M)$ in Theorem~\ref{thm:intersection-cohomology-ring}.

For equivariant cohomology, Theorem~\ref{thm:str-of-eq-cohomology} identifies $H_T^\ast(X)$ with the equivariant cohomology of the labelled graph associated with the characteristic pair. We next define equivariant intersection cohomology using finite-dimensional Borel approximations and study it by means of the associated spectral sequence and equivariant change of perversity. Finally, Theorem~\ref{thm:Stanley-Reisner-description} gives a combinatorial description of the middle-perversity equivariant intersection cohomology ring, and Corollary~\ref{cor:eq-intersection-cohomology-of-M} identifies it with $H_T^\ast(M)$.

\section{Preliminaries on intersection cohomology}\label{sec:some-remarks}

In this section, we fix the conventions for intersection cohomology used throughout Part~\ref{part-2} and recall the basic facts needed in the subsequent computations.
Throughout this part, all ordinary and intersection cohomology groups are taken with coefficients in $\mathbb Z$.
For intersection cohomology, we use Goresky--MacPherson (GM) perversities. A perversity specifies the allowable behavior of intersection chains near the singular strata of a stratified pseudomanifold.
More precisely, it assigns an integer to each possible codimension of a singular stratum, and these integers determine the allowability conditions in intersection homology.

A GM perversity $\bar p$ is a function on the integers $k\geq2$ satisfying
\[
\bar p(2)=0,
\qquad
\bar p(k)\leq \bar p(k+1)\leq \bar p(k)+1.
\]
Since the spaces considered in this paper are four-dimensional, only the values
$\bar p(2)$, $\bar p(3)$, and $\bar p(4)$ occur. 
Hence there are precisely four GM perversities, listed in Table~\ref{tab:GM-perversities-dim4}.
\begin{table}[htbp]
  \centering
  \begin{tabular}{c|ccc}
    \hline
    GM perversity & $\bar p(2)$ & $\bar p(3)$ & $\bar p(4)$ \\
    \hline
    $\bar 0$ & 0 & 0 & 0 \\
    $\bar m$ & 0 & 0 & 1 \\
    $\bar n$ & 0 & 1 & 1 \\
    $\bar t$ & 0 & 1 & 2 \\
    \hline
  \end{tabular}
  \caption{GM perversities in dimension four.}
  \label{tab:GM-perversities-dim4}
\end{table}
We call $\bar m$ and $\bar n$ the lower-middle and upper-middle perversities, respectively.
In this paper, middle-perversity intersection cohomology refers to intersection cohomology with respect to the lower-middle perversity $\bar m$.

For GM perversities, intersection cohomology is independent
of the choice of stratification, provided that there are
no codimension-one strata; see
\cite[Theorem~7.1.17]{Fri20}.
We use the stratification given by the orbit filtration
in Section~\ref{sec:locally-standard-T-pseudomanifold}.
We will frequently pass from intersection homology to intersection cohomology by the universal coefficient theorem \cite[Theorem~7.1.4]{Fri20}.
The arguments used in the computations below also show that all relevant intersection homology groups are free abelian.
Hence the $\operatorname{Ext}$ terms vanish, and we will use the resulting isomorphisms without further comment.

Recall that a stratified pseudomanifold is normal if all of its links are connected.
Since $X$ is a normal stratified pseudomanifold, top-perversity intersection cohomology agrees with ordinary cohomology (see \cite[Proposition~7.1.18]{Fri20}).
Thus, throughout the remainder of the paper, we identify
\[
I_{\bar t}H^\ast(X)
\cong
H^\ast(X).
\]

We use the orientation on $X$ fixed in Section~\ref{sec:setting}.
The Poincar\'e duality theorem for intersection homology \cite[Theorem~8.2.4]{Fri20} gives duality between complementary perversities.
In dimension four, the complementary pairs are $(\bar 0,\bar t)$ and $(\bar m,\bar n)$.
Consequently, the intersection cohomology groups for $\bar 0$ are determined by those for $\bar t$, while those for $\bar n$ are determined by those for $\bar m$.

If $M$ is a manifold, then, for every GM perversity $\bar p$, intersection cohomology agrees with ordinary cohomology (see \cite[Proposition~7.1.18]{Fri20}):
\[
I_{\bar p}H^\ast(M)\cong H^\ast(M).
\]

\section{Cohomology of a Thom block}
\label{sec:cohomology-of-Thom-block}

In this section, we compute the ordinary cohomology ring and the middle-perversity intersection cohomology groups of a Thom block.
We first determine its homotopy type in Proposition~\ref{prop:homotopy-type-of-Thom-block}, and then compute its intersection cohomology groups in Proposition~\ref{prop:intersection-cohomology-of-Thom-block}.

\begin{proposition}[homotopy type]\label{prop:homotopy-type-of-Thom-block}
The Thom space $\Th$ has the homotopy type
  \[
  \Th
  \simeq
  S^2\vee S^4.
  \]
\end{proposition}

\begin{proof}
For a space $Y$, let $Y_+:=Y\sqcup\{\ast\}$ denote the space obtained by adjoining a disjoint basepoint.
By \cite[Chapter~23, Section~5]{May99},
the Thom space of the trivial real $2$-plane bundle is
\[
\Th
\cong
Th(\C P^1\times\R^2)
\cong
(\C P^1)_+\wedge S^2
\cong
\Sigma^2((\C P^1)_+).
\]
Therefore, since $\C P^1\cong S^2$,
\[
\Sigma^2((\C P^1)_+)
\simeq
\Sigma^2\C P^1\vee S^2
\cong
S^4\vee S^2.
\]
Hence, we conclude
\[
\Th\simeq S^4\vee S^2.
\]
\end{proof}

It follows immediately from the previous proposition that we can determine the cohomology ring.

\begin{proposition}\label{prop:cohomology-ring-of-Thom-block}
  The cohomology ring of the Thom block is given by
  \[
  H^{\ast}(\Th)
\cong
\mathbb Z[u,b]/(u^2,ub,b^2),
\qquad
|u|=2,
\quad
|b|=4.
  \]
  As a module,
  \[
  H^i(\Th)\cong
  \begin{cases}
  \Z, & (i=0,2,4),\\
  0, & (\text{otherwise}),
\end{cases}
  \]
  and every product of positive-degree elements is zero.
\end{proposition}

\begin{remark}\label{rem:zero-top-Thom-block}
By Section~\ref{sec:some-remarks} and
Proposition~\ref{prop:cohomology-ring-of-Thom-block},
\[
I_{\bar 0}H^i(\Th)
\cong
I_{\bar t}H^i(\Th)
\cong
H^i(\Th)
\]
as abelian groups.
\end{remark}

\begin{remark}\label{rem:Thom-block-homotopy-warning}
Although $\Th\simeq S^2\vee S^4$, this homotopy equivalence cannot be used to compute the intersection homology of $\Th$.
Indeed, intersection homology is not a homotopy invariant in general.
Moreover, $S^2\vee S^4$ is not a $4$-dimensional stratified pseudomanifold, since its top-dimensional stratum is not dense.
\end{remark}

We next compute the middle-perversity intersection cohomology groups of the Thom block.

\begin{proposition}[intersection cohomology groups]\label{prop:intersection-cohomology-of-Thom-block}
For the lower middle perversity $\bar{m}$ and the upper middle perversity $\bar{n}$,
  \[
  I_{\bar{m}}H^{i}(\Th)
  \cong
  I_{\bar{n}}H^{i}(\Th)
  \cong
  \begin{cases}
    \Z, & (i=0,4),\\
    0, & (\text{otherwise}).
  \end{cases}
  \]
  Equivalently,
  \[
  I_{\bar{m}}H^{i}(\Th)
  \cong
  I_{\bar{n}}H^{i}(\Th)
  \cong
  H^i(S^4).
  \]
\end{proposition}

\begin{proof}
Let $X=Th(\mathbb{C}P^1\times\mathbb{C})$, equipped with the orbit filtration $X=X_{4}\supset X_{2}\supset X_{0}$. 
Let $v$ be its unique singular point, and let $L=\mathbb{C}P^{1}\times S^{1}$ be the link of $v$.

Choose an open neighborhood $V$ of $v$ and an open subset $U$ such that $X=U\cup V$, where $U$ is stratum-preserving homeomorphic to $\mathbb{C}P^1\times \mathring{D}^2$ (where $\mathring{D}^{2}$ is an open disk), $V=\mathring{c}(L)$ is an open cone on $L$, and $U\cap V$ is stratum-preserving homeomorphic to $L\times (0,1)$.

The induced codimension-two skeleta in $U$ and $U\cap V$ consist of smooth points. By stratification independence for GM perversities \cite[Theorem 7.1.17]{Fri20}, their intersection cohomology agrees with ordinary cohomology. 
Thus
\[
I_{\bar{p}}H^{i}(U)\cong H^{i}(\mathbb{C}P^{1}),\qquad
I_{\bar{p}}H^{i}(U\cap V)\cong H^{i}(L)\cong H^{i}(\mathbb{C}P^{1}\times S^{1}).
\]

For $\bar{p}=\bar{m}$ or $\bar{n}$, we have $\bar{p}(4)=1$. 
Since the formal dimension of the filtered space $L=\mathbb{C}P^{1}\times S^{1}$ is three, 
the cone formula \cite[Proposition 7.1.5]{Fri20}, applied to the induced filtration on $V$, gives
\begin{align*}
I_{\bar{p}}H^{i}(V)\cong
\begin{cases}
0, & (i>2), \\
{\rm Ext}(H_{1}(L),\mathbb{Z}), & (i=2), \\
H^{i}(L), & (i<2).
\end{cases}
\end{align*}
Because $H_{1}(L)\cong \mathbb{Z}$, we have $I_{\bar{p}}H^{2}(V)=0$. 
Moreover, the isomorphism $I_{\bar{p}}H^{1}(V)\to I_{\bar{p}}H^{1}(U\cap V)\cong H^{1}(L)$ is induced by the inclusion $U\cap V\hookrightarrow V$.

Since the inclusion $\mathbb CP^1\times S^1\hookrightarrow\mathbb CP^1\times \mathring{D}^2$ induces an isomorphism in degree two, the restriction
\[
I_{\bar p}H^2(U)\longrightarrow I_{\bar p}H^2(U\cap V)
\]
is also an isomorphism.

Since the restriction maps are isomorphisms in degrees one and two, the Mayer-Vietoris sequence for intersection cohomology \cite[Theorem 7.1.13]{Fri20} yields 
\begin{align*}
I_{\bar p}H^{i}(X)\cong
\begin{cases}
\mathbb Z, & (i=0,4), \\
0, & (\text{otherwise}).
\end{cases}
\end{align*}
This proves the proposition.
\end{proof}

\section{Ordinary cohomology ring}\label{sec:ordinary-cohomology-ring}

We now determine the ordinary cohomology ring of $X$ by combining the decomposition in Corollary~\ref{cor:decomposition-M-Thom} with the connected-sum formula in Theorem~\ref{thm:connected-sum-for-cohomology-ring} and the cohomology computation for a Thom block in Proposition~\ref{prop:cohomology-ring-of-Thom-block}.

\begin{theorem}\label{thm:ordinary-cohomology-ring}
There exist $u_1,\ldots,u_s\in H^2(X)$ such that
\[
H^{\ast}(X)
\cong
H^{\ast}(M)\oplus
\bigoplus_{j=1}^s \Z u_j
\]
as graded rings, where
$
u_i u_j=0
$
for all $i,j$, and
$
a u_j=0
$
for every $a\in H^{>0}(M)$.
\end{theorem}
\begin{proof}
This follows immediately from Theorem~\ref{thm:connected-sum-for-cohomology-ring} and
Proposition~\ref{prop:cohomology-ring-of-Thom-block}.
\end{proof}

\section{Intersection cohomology groups for all GM perversities}
\label{sec:intersection-cohomology-groups}

In this section, we determine the intersection cohomology groups of $X$ for all GM perversities.
We combine the connected-sum decomposition of Corollary~\ref{cor:decomposition-M-Thom} with Theorem~\ref{thm:connected-sum-for-intersection-cohomology} and the computation for a Thom block in Proposition~\ref{prop:intersection-cohomology-of-Thom-block}.
By Corollary~\ref{cor:decomposition-M-Thom}, we may write
\[
X \cong M \# B_1 \# \cdots \# B_s,
\]
where $M$ is a quasitoric manifold or $S^4$, and each $B_j$ is either $\Th$ or $\overline{\Th}$.
Although the connected sums are taken at regular fixed points,
these points have manifold neighborhoods.
Hence, by stratification independence for GM perversities,
Theorem~\ref{thm:connected-sum-for-intersection-cohomology}
applies.

We now state the intersection cohomology groups of $X$.

\begin{theorem}\label{thm:intersection-cohomology-group}
Let $X$ be a $4$-dimensional locally standard $T^2$-pseudomanifold satisfying Condition~\ref{cond:01-condition}.
Let $g$ be the number of vertices, $s$ the number of $0$-vertices, and $r$ the number of regular vertices, so that $r=g-s$.
Then
  \[
  I_{\bar{0}}H^i(X)
  \cong
  I_{\bar{t}}H^i(X)
  \cong
  H^i(X)
  \cong
  \begin{cases}
    \Z, &(i=0,4),\\
    \Z^{g-2}, &(i=2), \\
    0, &(\text{otherwise}),
  \end{cases}
  \]
and
  \[
  I_{\bar{m}}H^i(X)
  \cong
  I_{\bar{n}}H^i(X)
  \cong
  H^i(M)
  \cong
  \begin{cases}
    \Z, &(i=0,4),\\
    \Z^{r-2}=\Z^{g-s-2}, &(i=2), \\
    0, &(\text{otherwise}).
  \end{cases}
  \]
\end{theorem}

\begin{proof}
The assertion follows immediately from
Corollary~\ref{cor:decomposition-M-Thom},
Remark~\ref{rem:zero-top-Thom-block},
Proposition~\ref{prop:intersection-cohomology-of-Thom-block},
and Theorem~\ref{thm:connected-sum-for-intersection-cohomology}.
\end{proof}

The following corollary follows immediately from the theorem.

\begin{corollary}
Assume that $r \ge 3$, so that $M$ is a quasitoric manifold over an $r$-gon.
Then, for $i=0,1,2$,
\[
\rank H^{2i}(X)=h_i(P),\qquad
\rank (I_{\bar{m}}H^{2i}(X))
=
h_i(M/T),
\]
where $h(P)$ and $h(M/T)$ denote the $h$-vectors of
$P$ and $M/T$, respectively.
\end{corollary}

\section{Middle-perversity intersection cohomology ring}
\label{sec:intersection-cohomology-rings}

In general, intersection cohomology of a fixed perversity does not carry a ring structure induced by the intersection cup product, since the product may change the perversity.
In the present setting, however, the change of perversity from top to middle perversity allows us to define a ring structure on $I_{\bar m}H^\ast(X)$.
For perversities $\bar p\leq\bar q$, there is a natural change of perversity
$I_{\bar q}H^\ast(X)\to I_{\bar p}H^\ast(X)$ (see \cite[Proposition~7.1.7]{Fri20}).
Since $\bar m\leq\bar t$, and since we identify $I_{\bar t}H^\ast(X)\cong H^\ast(X)$
as in Section~\ref{sec:some-remarks}, we obtain a homomorphism
\[
\rho:
H^\ast(X)
\cong
I_{\bar t}H^\ast(X)
\to
I_{\bar m}H^\ast(X).
\]
We will show that $\rho$ is surjective and that its kernel is a homogeneous ideal of $H^\ast(X)$.
This allows us to equip $I_{\bar m}H^\ast(X)$ with the quotient ring structure induced from $H^\ast(X)$.
We then determine the kernel explicitly and show that the resulting intersection cohomology ring is isomorphic to $H^\ast(M)$, where $M$ is as in Corollary~\ref{cor:decomposition-M-Thom}.

\begin{remark}\label{rem:change-of-perversity}
For $\bar p\leq\bar q$, change of perversity is induced at the cochain level (see \cite[Proposition~7.1.7]{Fri20}).
Hence, by the standard naturality properties of cohomology, change of perversity is compatible with the long exact sequence of a pair, excision and the Mayer--Vietoris sequence.
We will use these compatibilities without further comment.
\end{remark}

We first determine the change of perversity $\rho$.

\begin{proposition}\label{prop:change-of-perversity}
The change of perversity
$\rho:H^\ast(X)\cong I_{\bar t}H^\ast(X)\to I_{\bar m}H^\ast(X)$
is surjective.
Moreover, if $X\cong M\# B_1\#\cdots\# B_s$
is the decomposition of Corollary~\ref{cor:decomposition-M-Thom}, and
$u_1,\ldots,u_s\in H^2(X)$ are the degree-$2$ classes corresponding to the Thom-block summands, then
\[
\ker\rho
=
\bigoplus_{j=1}^s\mathbb Z u_j
=
(u_1,\ldots,u_s).
\]
In particular, $\ker\rho$ is a homogeneous ideal of $H^\ast(X)$.
\end{proposition}

\begin{proof}
By Remark~\ref{rem:connected-sum-change-of-perversity}, change of perversity is compatible with the connected-sum decomposition. Hence the following diagram commutes:
\[
\begin{tikzcd}[ampersand replacement=\&]
H^2(X) \rar{\cong} \dar{\rho^2}
\&
H^2(M)\oplus\displaystyle\bigoplus_{j=1}^{s}H^2(B_j)
\dar{\rho^2}
\\
I_{\bar m}H^2(X) \rar{\cong}
\&
I_{\bar m}H^2(M)\oplus
\displaystyle\bigoplus_{j=1}^{s}I_{\bar m}H^2(B_j).
\end{tikzcd}
\]
Here the right vertical map is the direct sum of the corresponding changes of perversity on the connected-sum factors.
Since $M$ is a manifold, the change of perversity on the $M$-summand is the identity under the identification
$I_{\bar m}H^2(M)\cong H^2(M)$.
On the other hand, $I_{\bar m}H^2(B_j)=0$ for every $j$
by Proposition~\ref{prop:intersection-cohomology-of-Thom-block}.
Therefore, under the above identifications,
\[
\rho^2:
H^2(M)\oplus\bigoplus_{j=1}^s H^2(B_j)
\to H^2(M)
\]
is the projection onto the first factor. Hence
\[
\ker\rho^2
=
\bigoplus_{j=1}^s H^2(B_j)
=
\bigoplus_{j=1}^s\mathbb Z u_j.
\]
The change of perversity is an isomorphism in degrees $0$ and $4$ \cite[Corollary~8.1.26]{Fri20}, and all odd-degree groups vanish by Theorem~\ref{thm:intersection-cohomology-group}. Together with the description of $\rho^2$ above, this shows that $\rho$ is surjective and
\[
\ker\rho
=
\bigoplus_{j=1}^s\mathbb Z u_j.
\]
By Theorem~\ref{thm:connected-sum-for-cohomology-ring} and Proposition~\ref{prop:cohomology-ring-of-Thom-block}, we have
$u_i u_j=0$ for all $i,j$, and $a u_j=0$ for every
$a\in H^{>0}(M)$.
Hence
\[
\bigoplus_{j=1}^s\mathbb Z u_j
=
(u_1,\ldots,u_s),
\]
so $\ker\rho$ is a homogeneous ideal of $H^\ast(X)$.
\end{proof}

Proposition~\ref{prop:change-of-perversity} allows us to introduce the following ring structure.

\begin{definition}\label{def:ring-structure-of-intersection-cohomology}
By Proposition~\ref{prop:change-of-perversity}, the change of perversity induces an isomorphism of graded abelian groups
\[
H^\ast(X)/\ker\rho
\xrightarrow{\cong}
I_{\bar m}H^\ast(X).
\]
We equip $I_{\bar m}H^\ast(X)$ with the ring structure transported through this isomorphism from the quotient ring
$H^\ast(X)/\ker\rho$.
We call the resulting graded ring the \emph{middle-perversity intersection cohomology ring} of $X$.
\end{definition}

\begin{remark}
The ring structure in Definition~\ref{def:ring-structure-of-intersection-cohomology} is determined by the canonical change of perversity $\rho$ and therefore does not depend on the choice of the generators $u_1,\ldots,u_s$.
The generators are used only to describe $\ker\rho$ explicitly.
\end{remark}

We can now identify the resulting ring.

\begin{theorem}\label{thm:intersection-cohomology-ring}
Let $M$ be as in Corollary~\ref{cor:decomposition-M-Thom}.
Then there is an isomorphism of graded rings
\[
I_{\bar m}H^\ast(X)
\cong
H^\ast(M).
\]
\end{theorem}

\begin{proof}
By Proposition~\ref{prop:change-of-perversity} and Definition~\ref{def:ring-structure-of-intersection-cohomology},
\[
I_{\bar m}H^\ast(X)
\cong
H^\ast(X)/(u_1,\ldots,u_s).
\]
The cohomology-ring computation above gives
\[
H^\ast(X)
\cong
H^\ast(M)\oplus
\bigoplus_{j=1}^s\mathbb Z u_j,
\]
where all products involving a positive-degree element of $H^\ast(M)$ and some $u_j$, as well as all products $u_i u_j$, vanish.
Therefore quotienting by $(u_1,\ldots,u_s)$ gives
\[
H^\ast(X)/(u_1,\ldots,u_s)
\cong
H^\ast(M)
\]
as graded rings.
\end{proof}

\begin{remark}
For the orbit filtration of $X$, the only positive
codimensions are $2$ and $4$.
Since $\bar m(2)=\bar n(2)=0$ and
$\bar m(4)=\bar n(4)=1$,
the lower and upper middle perversities determine
the same intersection chain complex on $X$.

Hence
\[
I_{\bar m}H^\ast(X)\cong I_{\bar n}H^\ast(X).
\]
Under this identification, the changes of perversity
from $\bar t$ to $\bar m$ and $\bar n$ agree.
We equip both groups with the quotient ring structures
induced by these maps, as in
Definition~\ref{def:ring-structure-of-intersection-cohomology}.
The identification is then an isomorphism of graded rings.
\end{remark}

\subsection{An alternative description of the intersection cohomology product}
\label{subsec:alternative-intersection-product}

In this subsection, for perversities $\bar p\leq\bar q$, write
\[
\rho_{\bar q,\bar p}:
I_{\bar q}H^\ast(X)\to I_{\bar p}H^\ast(X)
\]
for the change of perversity.

The ring structure on $I_{\bar m}H^\ast(X)$ introduced in
Definition~\ref{def:ring-structure-of-intersection-cohomology}
is induced from the ordinary cup product $\smile$ on $H^\ast(X)$ through the quotient
$H^\ast(X)/\ker\rho$.
We now give an alternative description of its degree-$2$ product
using the intersection cup product.

\begin{construction}\label{cons:alternative-product}
The intersection cup product (see \cite[Section~7.3]{Fri20}) gives
\[
\smile_{\bar m,\bar m}:
I_{\bar m}H^2(X)\otimes I_{\bar m}H^2(X)
\to
I_{\bar 0}H^4(X).
\]
By \cite[Corollary~8.1.26]{Fri20}, change of perversity preserves the
fundamental cohomology class.
Hence
$\rho_{\bar m,\bar 0}:I_{\bar m}H^4(X)\to I_{\bar 0}H^4(X)$
is an isomorphism.
Composing $\smile_{\bar m,\bar m}$ with
$\rho_{\bar m,\bar 0}^{-1}$ defines a product
\[
\star:
I_{\bar m}H^2(X)\otimes I_{\bar m}H^2(X)
\xrightarrow{\smile_{\bar m,\bar m}}
I_{\bar 0}H^4(X)
\xrightarrow{\rho_{\bar{m},\bar{0}}^{-1}}
I_{\bar m}H^4(X).
\]
\end{construction}

We now compare this product with the quotient-ring product defined above.

\begin{proposition}
The product $\star$ defined in
Construction~\ref{cons:alternative-product}
agrees with the product on $I_{\bar m}H^\ast(X)$
defined in Definition~\ref{def:ring-structure-of-intersection-cohomology}.
\end{proposition}

\begin{proof}
Take $x,y\in I_{\bar m}H^2(X)$.
Since $\rho_{\bar t,\bar m}^2$ is surjective by Proposition~\ref{prop:change-of-perversity}, there exist
$a,b\in H^2(X)=I_{\bar t}H^2(X)$ such that
$\rho_{\bar t,\bar m}^2(a)=x$ and $\rho_{\bar t,\bar m}^2(b)=y$.
By Definition~\ref{def:ring-structure-of-intersection-cohomology},
the product of $x$ and $y$ is
$\rho_{\bar t,\bar m}^4(a\smile b)$.

By the naturality of the intersection cup product
\cite[Proposition~7.3.5]{Fri20}
and its compatibility with the ordinary cup product
\cite[Proposition~7.3.9]{Fri20},
the following diagram commutes:
\[
\begin{tikzcd}[ampersand replacement=\&]
H^2(X)\otimes H^2(X)
\ar[r,"\smile"]
\ar[d,"\rho_{\bar t,\bar m}^2\otimes\rho_{\bar t,\bar m}^2"']
\&
H^4(X)
\ar[d,"\rho_{\bar t,\bar 0}^4"]
\\
I_{\bar m}H^2(X)\otimes I_{\bar m}H^2(X)
\ar[r,"\smile_{\bar m,\bar m}"']
\&
I_{\bar 0}H^4(X).
\end{tikzcd}
\]
Hence
$x\smile_{\bar m,\bar m}y
=\rho_{\bar t,\bar 0}^4(a\smile b)$.
Since
$\rho_{\bar t,\bar 0}^4
=\rho_{\bar m,\bar 0}^4\circ\rho_{\bar t,\bar m}^4$,
we obtain
\[
(\rho_{\bar m,\bar 0}^4)^{-1}
(x\smile_{\bar m,\bar m}y)
=
\rho_{\bar t,\bar m}^4(a\smile b).
\]
The left-hand side is $x\star y$, while the right-hand side is the
product defined in
Definition~\ref{def:ring-structure-of-intersection-cohomology}.
Thus the two products agree.
\end{proof}

\section{Equivariant cohomology ring and the labelled graph}
\label{sec:equivariant-cohomology-labelled-graph}

In this section, we describe the equivariant cohomology ring of $X$ in terms of a labelled graph.
We continue to use the notation of Section~\ref{sec:setting}.
Thus, $P$ is a $g$-gon with edges $F_1,\ldots,F_g$ listed in cyclic order,
$v_i=F_i\cap F_{i+1}$, and $K_i=\lambda(F_i)$.
Indices are understood modulo $g$.
We use the canonical identification $\Hom(T,S^1)\cong H^2(BT)$.

The description developed in this section is motivated by the
GKM description of equivariant cohomology
\cite{GKM98,GZ01}.
In particular, we use the labelled-graph formulation developed by
Guillemin and Zara \cite{GZ01}.
In the present setting, however, $0$-vertices may occur, and the labels
on the two edges meeting at a $0$-vertex are linearly dependent.
Thus the labelled graph introduced below need not be a GKM graph in
the usual sense.

\subsection{The labelled graph associated with \texorpdfstring{$(P,\lambda)$}{(P, lambda)}}
\label{subsec:labelled-graph}
In this subsection, we define the labelled graph associated with $(P,\lambda)$ and its graph equivariant cohomology.
For each edge $F_i$, the characteristic functor $\lambda$ assigns a
circle subgroup $K_i:=\lambda(F_i)\subset T$.
Choose a primitive character
\[
\alpha_i\in\Hom(T,S^1)\cong H^2(BT)
\]
such that
\[
\ker\alpha_i=K_i.
\]
The character $\alpha_i$ is uniquely determined up to sign.

We now define the labelled graph associated with $(P,\lambda)$.

\begin{definition}[labelled graph]\label{def:labelled-graph}
The \emph{labelled graph} associated with $(P,\lambda)$ is the
$1$-skeleton $\Gamma$ of $P$ whose edge $F_i$ is labelled by
$\alpha_i$.
\end{definition}

We use these labels to define the corresponding graph equivariant cohomology.

\begin{definition}[graph equivariant cohomology]
\label{def:graph-eq-cohomology-of-Gamma}
The \emph{graph equivariant cohomology} of $\Gamma$ is
\[
H_T^\ast(\Gamma)
:=
\left\{
(f_1,\ldots,f_g)\in H^\ast(BT)^g
\ \middle|
f_i-f_{i-1}\in(\alpha_i)
\text{ for every }i
\right\},
\]
where $(\alpha_i)\subset H^\ast(BT)$ denotes the principal ideal
generated by $\alpha_i$.
With componentwise addition and multiplication,
$H_T^\ast(\Gamma)$ is a graded $H^\ast(BT)$-algebra, where
$H^\ast(BT)$ acts diagonally on $H^\ast(BT)^g$.
\end{definition}

\begin{remark}
The algebra $H_T^\ast(\Gamma)$ does not depend on the choices of signs
for the characters $\alpha_i$, since replacing $\alpha_i$ by
$-\alpha_i$ does not change the ideal $(\alpha_i)$.
\end{remark}

\begin{remark}
At a regular vertex $v_i$, the two incident labels
$\alpha_i$ and $\alpha_{i+1}$ form a basis of $H^2(BT)$, up to
the choices of signs.
At a $0$-vertex, on the other hand, $\alpha_{i+1}=\pm\alpha_i$,
and hence $(\alpha_i)=(\alpha_{i+1})$.
Thus $\Gamma$ need not satisfy the usual GKM independence
condition at every vertex.
\end{remark}

\subsection{Fixed-point description of equivariant cohomology}
\label{subsec:fixed-point-equivariant-cohomology}
In this subsection, we describe $H_T^\ast(X)$ by restriction to the fixed-point set.
Let $X^{T}=\{x_1,\ldots,x_g\}$ be the fixed-point set of the $T$-action on $X$ with $\pi(x_i)=v_i$.
We set
\[
X^{(1)}:=\{
x \in X \mid \dim T(x) \le 1
\}
=
\pi^{-1}(\partial P).
\]
For each $i$, put 
\[
Y_i:=\pi^{-1}(F_i).
\]
Then
\[
X^{(1)}=\bigcup_{i=1}^{g}Y_i,
\quad
Y_i \cap X^{T}=\{x_{i-1},x_i  \}.
\]

We first establish the collapse of the Borel spectral sequence.

\begin{lemma}\label{lem:Borel-spectral-sequence-collapse}
  The Serre spectral sequence of the Borel fibration
  \[
  X \to ET \times_T X \to BT
  \]
  degenerates at the $E_2$-page.
  In particular, the restriction to the fiber $\HT(X) \to H^{\ast}(X)$ is surjective.
\end{lemma}

\begin{proof}
Since $H^{\ast}(X)$ is a finitely generated free abelian group and
$H^{\mathrm{odd}}(X)=0$ by Theorem~\ref{thm:intersection-cohomology-group}, we have
\[
E_{2}^{a,b}
=
H^a(BT;H^{b}(X))
\cong
H^a(BT)\otimes H^b(X).
\]
Since both $H^{\mathrm{odd}}(BT)$ and $H^{\mathrm{odd}}(X)$ vanish,
the spectral sequence collapses at the $E_2$-page.
Hence $E_{\infty}^{0,b}=E_{2}^{0,b}\cong H^b(X)$.
Therefore, the restriction to the fiber $H_T^\ast(X)\to H^\ast(X)$ is surjective.
\end{proof}

We now apply the integral Chang--Skjelbred theorem to $X$.

\begin{proposition}[Chang--Skjelbred exact sequence]
  There is an exact sequence
  \[
  0 \to
  \HT(X)
  \xrightarrow{\iota^{\ast}}
  \HT(X^{T})
  \xrightarrow{\delta}
  H^{\ast+1}_{T}(X^{(1)},X^{T}),
  \]
  where $\iota:X^{T} \to X$ is the inclusion.
  In particular, 
  $\HT(X)
  \xrightarrow{\iota^{\ast}}
  \HT(X^{T})$
  is injective.
\end{proposition}

\begin{proof}
By the canonical-model identification, we may regard $X$ as $X(P,\lambda)$.
The face decomposition of $P$ induces a finite $T$-CW structure on $X$, since over the relative interior of each face $F$ the corresponding stratum is $T$-equivariantly homeomorphic to $(T/\lambda(F))\times\mathring F$.
By Lemma~\ref{lem:Borel-spectral-sequence-collapse}, the restriction $H_T^\ast(X)\to H^\ast(X)$ is surjective.
Moreover, every isotropy subgroup of $X$ is a subtorus of $T$, and hence is connected. Thus the hypotheses of \cite[Theorem~1.1]{FP07} are satisfied, and the integral Chang--Skjelbred sequence is exact. Hence
\[
0\to
H_T^\ast(X)
\to
H_T^\ast(X^T)
\to
H_T^{\ast+1}(X^{(1)},X^T)
\]
is exact.
\end{proof}

Since every $x_i$ is a fixed point, we have canonical identifications
\[
\HT(X^{T})
\cong
\bigoplus_{i=1}^{g} \HT(x_i)
\cong
H^{\ast}(BT)^g.
\]
Under these identifications, the homomorphism
$\iota^\ast:H_T^\ast(X)\to H_T^\ast(X^T)$
is given by restriction to the fixed points:
\[
\iota^\ast(f)
=
(f|_{x_1},\ldots,f|_{x_g}).
\]

We next determine the restriction condition contributed by each edge.

\begin{lemma}[edge restriction]
\label{lem:edge-restriction}
For each edge $F_i$, the image of the restriction homomorphism
\[
H_T^\ast(Y_i)
\to
H_T^\ast(x_{i-1})\oplus H_T^\ast(x_i)
\cong
H^\ast(BT)^2
\]
is
\[
\left\{
(f_{i-1},f_i)\in H^\ast(BT)^2
\ \middle|
f_i-f_{i-1}\in(\alpha_i)
\right\}.
\]
\end{lemma}

\begin{proof}
This follows immediately from the Mayer--Vietoris sequence.
\end{proof}

We next decompose the relative equivariant cohomology of the $1$-skeleton into the contributions from the individual edges.

\begin{lemma}
There is a natural isomorphism
\[
H_T^\ast(X^{(1)},X^T)
\cong
\bigoplus_{i=1}^g
H_T^\ast(Y_i,\{x_{i-1},x_i\}).
\]
\end{lemma}

\begin{proof}
If $i\neq j$, then $Y_i\cap Y_j\subset X^T$.
Passing to the Borel constructions and collapsing $(X^T)_T$ to a point,
we obtain a natural homeomorphism
\[
(X^{(1)})_T/(X^T)_T
\cong
\bigvee_{i=1}^g
(Y_i)_T/(\{x_{i-1},x_i\})_T.
\]
Since relative cohomology is identified with the reduced cohomology of
the corresponding quotient, and reduced cohomology takes a finite
wedge to a direct sum, the asserted isomorphism follows.
\end{proof}

We now determine the image of the fixed-point restriction map.

\begin{theorem}[structure of equivariant cohomology]
\label{thm:str-of-eq-cohomology}
The fixed-point restriction map induces an isomorphism of graded
$H^\ast(BT)$-algebras
\[
H_T^\ast(X)
\cong
H_T^\ast(\Gamma)
=
\left\{
(f_1,\ldots,f_g)\in H^\ast(BT)^g
\ \middle|
f_i-f_{i-1}\in(\alpha_i)
\text{ for every }i
\right\}.
\]
\end{theorem}

\begin{proof}
This follows immediately from the Chang--Skjelbred exact
sequence and the preceding two lemmas.
\end{proof}

We next describe the $0$-vertices in terms of degree-$2$ equivariant cohomology classes supported at single fixed points.

\begin{proposition}
\label{prop:zero-vertex-class}
A vertex $v_i$ is a $0$-vertex if and only if there exists
$\zeta_i\in H_T^2(X)$ such that
\[
\zeta_i|_{x_i}=\alpha_i,
\qquad
\zeta_i|_{x_j}=0
\quad\text{for every }j\neq i.
\]
Under the fixed-point restriction isomorphism of
Theorem~\ref{thm:str-of-eq-cohomology}, this condition is written as
\[
\zeta_i=(0,\ldots,0,\alpha_i,0,\ldots,0).
\]
\end{proposition}

\begin{proof}
Suppose first that $v_i$ is a $0$-vertex.
Then $(\alpha_i)=(\alpha_{i+1})$.
Hence the tuple
\[
(0,\ldots,0,\alpha_i,0,\ldots,0)
\]
satisfies all the congruence
relations in Theorem~\ref{thm:str-of-eq-cohomology}.
Thus it defines a class $\zeta_i\in H_T^2(X)$.

Conversely, suppose that such a class $\zeta_i$ exists.
Applying the congruence relation corresponding to the edge
$F_{i+1}$ gives
$\alpha_i\in(\alpha_{i+1})$.
Since $\alpha_i$ and $\alpha_{i+1}$ are primitive characters, this
implies $\alpha_i=\pm\alpha_{i+1}$.
Thus $K_i=K_{i+1}$, and hence $v_i$ is a $0$-vertex.
\end{proof}

We conclude this section by describing how the ordinary cohomology ring is recovered from the equivariant cohomology ring.

\begin{corollary}
There is an isomorphism of graded rings
\[
H^\ast(X)
\cong
H_T^\ast(X)\otimes_{H^\ast(BT)}\mathbb Z,
\]
where $\mathbb Z$ is regarded as an $H^\ast(BT)$-module via the
augmentation $H^\ast(BT)\to\mathbb Z$.
\end{corollary}

\begin{proof}
By Lemma~\ref{lem:Borel-spectral-sequence-collapse}, the restriction
$H_T^\ast(X)\to H^\ast(X)$ is surjective.
Since $H^\ast(X)$ is a finitely generated free abelian group,
the Leray--Hirsch theorem gives an isomorphism of $H^\ast(BT)$-modules
\[
H_T^\ast(X)\cong H^\ast(BT)\otimes H^\ast(X).
\]
Regarding $H^\ast(X)$ as an $H^\ast(BT)$-algebra via the augmentation,
the restriction homomorphism factors through
\[
H_T^\ast(X)\otimes_{H^\ast(BT)}\mathbb Z
\longrightarrow
H^\ast(X).
\]
By the Leray--Hirsch isomorphism above, this homomorphism is an
isomorphism of graded abelian groups.
Since it is induced by the restriction homomorphism, it is an
isomorphism of graded rings.
\end{proof}

\section{Equivariant intersection cohomology}
\label{sec:equivariant-intersection-cohomology}

In this section, we define equivariant intersection cohomology and study its group and $H^\ast(BT)$-module structures, as well as the quotient ring structure in middle perversity.

\subsection{Definition of equivariant intersection cohomology}
\label{subsec:def-equivariant-intersection-cohomology}
In this subsection, we define equivariant intersection cohomology using finite-dimensional Borel approximations.
Let $T=T^m$. For each $N\geq0$, consider the finite-dimensional
approximations
\[
E_NT=(S^{2N+1})^m,
\qquad
B_NT=(\mathbb CP^N)^m,
\]
and set $X_{T,N}:=E_NT\times_T X$.
The stratification of $X$ induces a natural stratification on
$X_{T,N}$, and hence $I_{\bar p}H^\ast(X_{T,N})$ is defined for every
GM perversity $\bar p$.

For $N'\geq N$, the natural $T$-equivariant inclusion
$E_NT\hookrightarrow E_{N'}T$ induces a stratification-preserving
inclusion
\[
j_{N,N'}:X_{T,N}\to X_{T,N'}.
\]
These inclusions induce homomorphisms
\[
j_{N,N'}^\ast:
I_{\bar p}H^\ast(X_{T,N'})
\to
I_{\bar p}H^\ast(X_{T,N}),
\]
and hence the groups $I_{\bar p}H^\ast(X_{T,N})$ form an inverse
system.
We follow the finite-dimensional Borel approximation approach to
equivariant intersection cohomology used in \cite{Kir88}.

\begin{definition}
The \emph{$T$-equivariant intersection cohomology} of $X$ with
perversity $\bar p$ is defined by
\[
I_{\bar p}H_T^\ast(X)
:=
\varprojlim_N I_{\bar p}H^\ast(X_{T,N}).
\]
\end{definition}

\begin{remark}
For the top perversity $\bar t$, we identify
\[
I_{\bar t}H_T^\ast(X)
\cong
H_T^\ast(X).
\]
Indeed, for each $N$, since $X_{T,N}$ is a normal stratified
pseudomanifold,
\[
I_{\bar t}H^\ast(X_{T,N})
\cong
H^\ast(X_{T,N})
\]
by \cite[Proposition~7.1.18]{Fri20}.
These identifications are compatible with the transition maps.
Together with degree-wise stabilization and the usual finite-dimensional
approximation of Borel equivariant cohomology, this gives the asserted
identification.
\end{remark}

For the spectral-sequence argument below, we impose the following
condition.

\begin{condition}
\label{cond:odd-intersection-cohomology-vanish}
Let $X$ be a normal locally standard $T^m$-pseudomanifold and let
$\bar p$ be a GM perversity.
Suppose that $I_{\bar p}H^{\mathrm{odd}}(X)=0$ and that
$I_{\bar p}H^\ast(X)$ is a finitely generated free abelian group.
\end{condition}

Throughout the following spectral-sequence argument, we assume
Condition~\ref{cond:odd-intersection-cohomology-vanish}.
For the $4$-dimensional locally standard $T^2$-pseudomanifolds
considered in this paper, this condition is satisfied by
Theorem~\ref{thm:intersection-cohomology-group}.

\begin{remark}
  Under Condition~\ref{cond:odd-intersection-cohomology-vanish}, for every degree $k$, the inverse system 
\[
\{I_{\bar p}H^k(X_{T,N})\}_N
\]
stabilizes; see Lemma~\ref{lem:degree-wise-stabilization}. Hence, for every sufficiently large $N$, we have $I_{\bar p}H_T^k(X)\cong I_{\bar p}H^k(X_{T,N})$.
\end{remark}

\subsection{The equivariant intersection cohomology spectral sequence}
\label{subsec:equivariant-IH-spectral-sequence}
In this subsection, we study the spectral sequence associated with a filtration of the finite-dimensional Borel approximation.
Let $\pi_N:X_{T,N}\to B_NT$
be the projection of the finite-dimensional Borel approximation.
Fix a finite triangulation of $B_NT$, and let $S_a$ denote its $a$-skeleton.
Set
\[
F_aX_{T,N}:=\pi_N^{-1}(S_a).
\]
Then
\[
\emptyset=F_{-1}X_{T,N}
\subset
F_0X_{T,N}
\subset\cdots\subset
F_{D_N}X_{T,N}=X_{T,N},
\]
where $D_N=\dim B_NT$.
We consider the intersection cohomology spectral sequence associated with this filtration.
Its $E_1$-term is
\[
E_{1,\bar p}^{a,b}
\cong
I_{\bar p}H^{a+b}
(F_aX_{T,N},F_{a-1}X_{T,N}),
\]
and the spectral sequence converges to
$I_{\bar p}H^{a+b}(X_{T,N})$.
We now compute the $E_1$- and $E_2$-terms
and show that the spectral sequence collapses at $E_2$.

\begin{lemma}
\label{lem:E_2-term}
There are isomorphisms
\begin{align*}
E_{1,\bar p}^{a,b}
&\cong
C^a\bigl(B_NT;I_{\bar p}H^b(X)\bigr),\\
E_{2,\bar p}^{a,b}
&\cong
H^a(B_NT)\otimes I_{\bar p}H^b(X).
\end{align*}
The spectral sequence collapses at the $E_2$-page:
\[
E_{2,\bar p}=E_{\infty,\bar p}.
\]
\end{lemma}

\begin{proof}
Let $\Delta^a$ be an $a$-simplex of $B_NT$.
Since $\Delta^a$ is contractible and
$\pi_N:X_{T,N}\to B_NT$ is a bundle with fiber $X$,
its restriction over $\Delta^a$ is trivial.
By excision, the relative group
$E_{1,\bar p}^{a,b}\cong
I_{\bar p}H^{a+b}(F_aX_{T,N},F_{a-1}X_{T,N})$
is the direct sum of the groups
\[
I_{\bar p}H^{a+b}
(\Delta^a\times X,\partial\Delta^a\times X).
\]
By the K\"unneth formula and the freeness assumption in
Condition~\ref{cond:odd-intersection-cohomology-vanish},
each group is isomorphic to
\[
H^a(\Delta^a,\partial\Delta^a)
\otimes I_{\bar p}H^b(X).
\]
Since $B_NT$ is simply connected, the local coefficient
system is constant.
Taking the direct sum over all $a$-simplices gives
\[
E_{1,\bar p}^{a,b}
\cong
C^a\bigl(B_NT;I_{\bar p}H^b(X)\bigr).
\]
Under this identification, the differential $d_1$
is the usual coboundary.
Hence
\[
E_{2,\bar p}^{a,b}\cong H^a(B_NT;I_{\bar p}H^b(X)).
\]
Since the integral homology of $B_NT$ is free abelian and
$I_{\bar p}H^b(X)$ is free abelian by
Condition~\ref{cond:odd-intersection-cohomology-vanish},
the universal coefficient theorem gives
\[
H^a(B_NT;I_{\bar p}H^b(X))
\cong
H^a(B_NT)\otimes I_{\bar p}H^b(X).
\]
Since $H^{\mathrm{odd}}(B_NT)=0$,
Condition~\ref{cond:odd-intersection-cohomology-vanish}
implies that $E_{2,\bar p}^{a,b}=0$
unless both $a$ and $b$ are even.
The differential $d_r$ has bidegree $(r,1-r)$,
so it vanishes for every $r\geq2$.
\end{proof}

We next describe the filtration on the abutment of the spectral sequence.

\begin{lemma}[convergence of the spectral sequence]
\label{lem:convergence-equivariant-IH}
For each fixed degree $k$, there is a finite decreasing filtration
\[
I_{\bar p}H^k(X_{T,N})
=
\mathcal F_{\bar p}^{0}
\supset
\mathcal F_{\bar p}^{1}
\supset\cdots\supset
\mathcal F_{\bar p}^{D_N}
\supset
\mathcal F_{\bar p}^{D_N+1}=0,
\]
where
\[
\mathcal F_{\bar p}^{a}
:=
\operatorname{Im}
\left(
I_{\bar p}H^k(X_{T,N},F_{a-1}X_{T,N})
\to
I_{\bar p}H^k(X_{T,N})
\right).
\]
Moreover,
\[
\mathcal F_{\bar p}^{a}/\mathcal F_{\bar p}^{a+1}
\cong
E_{\infty,\bar p}^{a,k-a}.
\]
\end{lemma}

\begin{proof}
This follows from the convergence theorem for the cohomological
spectral sequence associated with a finite filtration
\cite[Theorem~2.6]{McC01}.
\end{proof}

\begin{remark}
By the long exact sequence of the pair
$(X_{T,N},F_{a-1}X_{T,N})$, one can equivalently write
\[
\mathcal F_{\bar p}^{a}
=
\ker
\left(
I_{\bar p}H^k(X_{T,N})
\to
I_{\bar p}H^k(F_{a-1}X_{T,N})
\right).
\]
\end{remark}

We next determine the abelian group structure of $I_{\bar p}H^\ast(X_{T,N})$.

\begin{lemma}[group structure]
\label{lem:group-structure-1}
For each $k$, there is an isomorphism of abelian groups
\[
I_{\bar p}H^k(X_{T,N})
\cong
\bigoplus_{a+b=k}E_{2,\bar p}^{a,b}
\cong
\bigoplus_{a+b=k}
H^a(B_NT)\otimes I_{\bar p}H^b(X).
\]
\end{lemma}

\begin{proof}
By the convergence of the spectral sequence and its collapse at the
$E_2$-page, we have
\[
\mathcal F_{\bar p}^{a}/\mathcal F_{\bar p}^{a+1}
\cong
E_{\infty,\bar p}^{a,k-a}
\cong
E_{2,\bar p}^{a,k-a}
\cong
H^a(B_NT)\otimes I_{\bar p}H^{k-a}(X).
\]
By Condition~\ref{cond:odd-intersection-cohomology-vanish}, the group
$I_{\bar p}H^{k-a}(X)$ is free abelian, and hence so is
$H^a(B_NT)\otimes I_{\bar p}H^{k-a}(X)$.
Therefore, the short exact sequence
\[
0\to
\mathcal F_{\bar p}^{a+1}
\to
\mathcal F_{\bar p}^{a}
\to
\mathcal F_{\bar p}^{a}/\mathcal F_{\bar p}^{a+1}
\to0
\]
splits.
Starting with $\mathcal F_{\bar p}^{D_N+1}=0$,
we split these sequences successively.
Since $\mathcal F_{\bar p}^{0}=I_{\bar p}H^k(X_{T,N})$,
we obtain
\[
I_{\bar p}H^k(X_{T,N})
\cong
\bigoplus_{a+b=k}E_{2,\bar p}^{a,b}
\cong
\bigoplus_{a+b=k}
H^a(B_NT)\otimes I_{\bar p}H^b(X).
\]
\end{proof}

\begin{remark}
The isomorphism in Lemma~\ref{lem:group-structure-1}
depends on the choice of splittings and is not canonical.
\end{remark}

We next define the $H^\ast(B_NT)$-module structure on $I_{\bar p}H^\ast(X_{T,N})$.

\begin{lemma}[module structure]
\label{lem:module-structure}
The abelian group $I_{\bar p}H^\ast(X_{T,N})$ has a natural
$H^\ast(B_NT)$-module structure.
\end{lemma}

\begin{proof}
Since $X_{T,N}$ is a normal pseudomanifold,
$I_{\bar t}H^\ast(X_{T,N})\cong H^\ast(X_{T,N})$.
Moreover, by \cite[Corollary~7.2.12]{Fri20},
$(\bar p,\bar t;\bar p)$ is agreeable, so the intersection cup product
\[
\smile_{\bar p,\bar t}:
I_{\bar p}H^i(X_{T,N})
\otimes
I_{\bar t}H^j(X_{T,N})
\to
I_{\bar p}H^{i+j}(X_{T,N})
\]
is defined.
For $a\in H^j(B_NT)$ and
$x\in I_{\bar p}H^i(X_{T,N})$, define
\[
a\cdot x
:=
x\smile_{\bar p,\bar t}\pi_N^\ast(a),
\]
where
$\pi_N^\ast(a)\in H^j(X_{T,N})
\cong I_{\bar t}H^j(X_{T,N})$.
The unitality and associativity of this action follow
from those of the intersection cup product
\cite[Propositions~7.3.21,~7.3.34]{Fri20}.
\end{proof}

\subsection{Structure of equivariant intersection cohomology}
\label{subsec:structure-equivariant-intersection-cohomology}

We first show that the inverse system stabilizes in each degree.

\begin{lemma}[degree-wise stabilization]
\label{lem:degree-wise-stabilization}
For every $k$, there exists $N_0$ such that, for all
$N'\geq N\geq N_0$, the natural map
\[
I_{\bar p}H^k(X_{T,N'})
\to
I_{\bar p}H^k(X_{T,N})
\]
is an isomorphism.
Consequently, for every sufficiently large $N$,
\[
I_{\bar p}H_T^k(X)
\cong
I_{\bar p}H^k(X_{T,N}).
\]
\end{lemma}

\begin{proof}
Fix $k$, and let $N'\geq N$. The natural inclusion
$\iota_{N,N'}:B_NT\to B_{N'}T$ induces
$j_{N,N'}:X_{T,N}\to X_{T,N'}$.
We show that, for $N$ sufficiently large, the induced homomorphism
$j_{N,N'}^\ast:I_{\bar p}H^k(X_{T,N'})\to I_{\bar p}H^k(X_{T,N})$
is an isomorphism.

Consider the spectral sequences associated with the filtrations of
$X_{T,N'}$ and $X_{T,N}$ introduced above. By
Lemma~\ref{lem:E_2-term}, their $E_2$-terms are
$H^a(B_{N'}T)\otimes I_{\bar p}H^b(X)$ and
$H^a(B_NT)\otimes I_{\bar p}H^b(X)$, respectively.
The map $j_{N,N'}^\ast$ induces a morphism between these spectral
sequences, and on the $E_2$-page it is given by
$\iota_{N,N'}^\ast\otimes 1$.
Since $B_NT=(\mathbb CP^N)^m$, the homomorphism
$\iota_{N,N'}^\ast:H^a(B_{N'}T)\to H^a(B_NT)$ is an isomorphism for
$a\leq 2N$. Choose $N$ so that $2N\geq k$. 
Write $E_{r,\bar p}^{a,b}(N)$ for the spectral sequence associated with $X_{T,N}$.
Then, for every $a+b=k$, the induced map
\[
E_{2,\bar p}^{a,b}(N')
\to
E_{2,\bar p}^{a,b}(N)
\]
is an isomorphism.
Both spectral sequences collapse at the $E_2$-page. Hence the induced
map on the corresponding $E_\infty$-terms is an isomorphism for every
$a+b=k$. 
By Lemma~\ref{lem:convergence-equivariant-IH}, these $E_\infty$-terms are the successive quotients of finite filtrations of $I_{\bar p}H^k(X_{T,N'})$ and $I_{\bar p}H^k(X_{T,N})$.
It follows by induction on the filtrations that
\[
j_{N,N'}^\ast:
I_{\bar p}H^k(X_{T,N'})
\xrightarrow{\cong}
I_{\bar p}H^k(X_{T,N}).
\]
Thus the inverse system stabilizes in degree $k$. Therefore, for every
sufficiently large $N$,
\[
I_{\bar p}H_T^k(X)
\cong
I_{\bar p}H^k(X_{T,N}).
\]
\end{proof}

We next define the $H^\ast(BT)$-module structure on equivariant intersection cohomology.

\begin{lemma}\label{lem:eq-IH-module-structure}
  The equivariant intersection cohomology
  $I_{\bar{p}}H^{\ast}_{T}(X)$
  has an $H^{\ast}(BT)$-module structure.
  The action of $H^{\ast}(BT)$ is induced by the $H^{\ast}(B_NT)$-action in Lemma~\ref{lem:module-structure}.
\end{lemma}

\begin{proof}
Define the $H^\ast(BT)$-action as follows.
Let $a\in H^j(BT)$ and
$x=(x_N)_N\in I_{\bar p}H_T^i(X)$.
For each $N$, let $\iota_N:B_NT\to BT$ be the natural inclusion and put
$a_N:=\iota_N^\ast(a)\in H^j(B_NT)$.
By the definition of equivariant intersection cohomology,
$x$ is a compatible family with
$x_N\in I_{\bar p}H^i(X_{T,N})$.
Using the $H^\ast(B_NT)$-module structure defined in
Lemma~\ref{lem:module-structure}, set
$(a\cdot x)_N:=a_N\cdot x_N \in I_{\bar p}H^{i+j}(X_{T,N})$.
We show that the family $(a_N\cdot x_N)_N$ is compatible with the
transition maps.

Let $N'\geq N$, and denote the natural inclusions by
\[
\iota_{N,N'}:B_NT\to B_{N'}T,
\qquad
j_{N,N'}:X_{T,N}\to X_{T,N'}.
\]
Since $\iota_N=\iota_{N'}\circ\iota_{N,N'}$, we have
$\iota_{N,N'}^\ast(a_{N'})=a_N$.
Moreover, since $x=(x_N)_N$ belongs to the inverse limit,
$j_{N,N'}^\ast(x_{N'})=x_N$.
The projections satisfy
\[
\pi_{N'}\circ j_{N,N'}
=
\iota_{N,N'}\circ\pi_N,
\]
and hence
\[
j_{N,N'}^\ast\pi_{N'}^\ast(a_{N'})
=
\pi_N^\ast(a_N).
\]
By the definition of the action in
Lemma~\ref{lem:module-structure} and the naturality of the
intersection cup product
\cite[Proposition~7.3.5]{Fri20}, we obtain
\[
\begin{aligned}
j_{N,N'}^\ast(a_{N'}\cdot x_{N'})
&=
j_{N,N'}^\ast
\left(
x_{N'}\smile_{\bar p,\bar t}\pi_{N'}^\ast(a_{N'})
\right)\\
&=
j_{N,N'}^\ast(x_{N'})
\smile_{\bar p,\bar t}
j_{N,N'}^\ast\pi_{N'}^\ast(a_{N'})\\
&=
x_N\smile_{\bar p,\bar t}\pi_N^\ast(a_N)\\
&=
a_N\cdot x_N.
\end{aligned}
\]
Thus $(a_N\cdot x_N)_N$ is a compatible family and defines an element
of $I_{\bar p}H_T^{i+j}(X)$.
The unitality and associativity of this action
follow from those at each finite stage.
\end{proof}

Combining degree-wise stabilization with the finite-stage computation, we obtain the abelian group structure of equivariant intersection cohomology.

\begin{lemma}[group structure of equivariant intersection cohomology]
\label{lem:eq-IH-group-structure}
For each $k$, there is an isomorphism of abelian groups
\[
I_{\bar p}H_T^k(X)
\cong
\bigoplus_{a+b=k}
H^a(BT)\otimes I_{\bar p}H^b(X).
\]
\end{lemma}

\begin{proof}
Fix $k$ and take $N$ sufficiently large with $2N\geq k$.
By Lemmas~\ref{lem:degree-wise-stabilization}
and~\ref{lem:group-structure-1}, we have
\[
I_{\bar p}H_T^k(X)
\cong
\bigoplus_{a+b=k}
H^a(B_NT)\otimes I_{\bar p}H^b(X).
\]
The inclusion $B_NT\to BT$ induces isomorphisms
$H^a(BT)\cong H^a(B_NT)$ for $a\leq2N$.
Since $a+b=k$ implies $a\leq k$, the assertion follows.
\end{proof}

\subsection{Change of perversity}\label{subsec:equivariant-change-of-perversity}
In this subsection, we study change of perversity from $\bar t$ to $\bar m$.
For each $N$, change of perversity induces a homomorphism
\[
\rho_N:
I_{\bar t}H^\ast(X_{T,N})
\to
I_{\bar m}H^\ast(X_{T,N}).
\]
We denote the nonequivariant change of perversity on $X$ by
\[
\rho:
I_{\bar t}H^\ast(X)
\to
I_{\bar m}H^\ast(X).
\]
We first describe the induced map on the $E_2$-term.

\begin{lemma}[change of perversity on the $E_2$-term]
\label{lem:change-of-perversity-on-E_2-term}
Under the identification
\[
E_{2,\bar p}^{a,b}
\cong
H^a(B_NT)\otimes I_{\bar p}H^b(X),
\]
the homomorphism induced by $\rho_N$ on the $E_2$-term is
\[
1\otimes\rho^b:
H^a(B_NT)\otimes I_{\bar t}H^b(X)
\to
H^a(B_NT)\otimes I_{\bar m}H^b(X).
\]
\end{lemma}
\begin{proof}
At the cochain level, we consider the map
$\rho_N^C:I_{\bar t}C^\ast(X_{T,N})
\to I_{\bar m}C^\ast(X_{T,N})$
inducing $\rho_N$
(see \cite[Proposition~7.1.7]{Fri20}).
For every $a$, this map restricts to the relative cochain
complexes for $(X_{T,N},F_{a-1}X_{T,N})$.
Thus it preserves the filtration and induces a morphism
of the associated spectral sequences:
\[
\rho_r^{a,b}:
E_{r,\bar t}^{a,b}
\to
E_{r,\bar m}^{a,b},
\qquad r\geq1.
\]

Let $\Delta^a$ be an $a$-simplex of $B_NT$.
Since $\pi_N^{-1}(\Delta^a)\cong\Delta^a\times X$, the relative
K\"unneth isomorphism gives
\[
I_{\bar p}H^{a+b}
(\Delta^a\times X,\partial\Delta^a\times X)
\cong
H^a(\Delta^a,\partial\Delta^a)
\otimes I_{\bar p}H^b(X).
\]
By the naturality of this isomorphism with respect to
change of perversity in the $X$-factor, the induced
homomorphism is $1\otimes\rho^b$.
Taking the direct sum over all $a$-simplices and using
Lemma~\ref{lem:E_2-term}, we obtain the map
\[
C^a(B_NT;\rho^b):
C^a(B_NT;I_{\bar t}H^b(X))
\to
C^a(B_NT;I_{\bar m}H^b(X)),
\]
where $C^a(B_NT;\rho^b)$ denotes the homomorphism obtained by applying
$\rho^b$ to the coefficient group on each $a$-simplex.

Since $d_1$ is the usual coboundary, the induced map
on the $E_2$-term is $H^a(B_NT;\rho^b)$.
The identifications in Lemma~\ref{lem:E_2-term} are natural
with respect to homomorphisms of the coefficient groups.
Hence this map agrees with $1\otimes\rho^b$, as claimed.
\end{proof}

The homomorphisms $\rho_N$ are compatible with the transition maps of
the finite-dimensional Borel approximations and therefore induce a homomorphism
\[
\rho_T:
I_{\bar t}H_T^\ast(X)
\to
I_{\bar m}H_T^\ast(X).
\]
We next prove that the equivariant change of perversity $\rho_T$ is surjective.

\begin{lemma}[surjectivity of equivariant change of perversity]
\label{lem:surjectivity-of-equivariant-change-of-perversity}
Assume that
Condition~\ref{cond:odd-intersection-cohomology-vanish}
holds for both $\bar t$ and $\bar m$, and that the
nonequivariant change of perversity $\rho$ is surjective.
Then the homomorphism
\[
\rho_T:
I_{\bar t}H_T^\ast(X)
\to
I_{\bar m}H_T^\ast(X)
\]
is surjective.
\end{lemma}

\begin{proof}
Fix $k$ and take $N$ sufficiently large.
By Lemma~\ref{lem:degree-wise-stabilization},
$\rho_T^k$ corresponds to $\rho_N^k$ under the stabilization
isomorphisms.
Thus it suffices to prove that $\rho_N^k$ is surjective.

By Lemma~\ref{lem:change-of-perversity-on-E_2-term}, the homomorphism
induced on the $E_2$-term is
\[
1\otimes\rho^b:
H^a(B_NT)\otimes I_{\bar t}H^b(X)
\to
H^a(B_NT)\otimes I_{\bar m}H^b(X).
\]
The nonequivariant change of perversity
$\rho^b:I_{\bar t}H^b(X)\to I_{\bar m}H^b(X)$
is surjective by assumption.
Hence $1\otimes\rho^b$ is also surjective.
The morphism of spectral sequences in
Lemma~\ref{lem:change-of-perversity-on-E_2-term}
induces a homomorphism
\[
\rho_\infty^{a,b}:
E_{\infty,\bar t}^{a,b}
\to
E_{\infty,\bar m}^{a,b}.
\]
Since both spectral sequences collapse at the $E_2$-page,
this map is identified with $1\otimes\rho^b$
and is therefore surjective.

We now use the finite filtrations from
Lemma~\ref{lem:convergence-equivariant-IH},
whose successive quotients satisfy
\[
\mathcal F_{\bar p}^{a}/\mathcal F_{\bar p}^{a+1}
\cong
E_{\infty,\bar p}^{a,k-a}.
\]
The homomorphism $\rho_N^k$ preserves these filtrations
and induces $\rho_\infty^{a,k-a}$ on the quotients.
These induced maps are therefore surjective.
We apply descending induction on $a$, starting with
$\mathcal F_{\bar p}^{D_N+1}=0$.
The short exact sequences of the filtrations show at each step
that
$\rho_N^k:\mathcal F_{\bar t}^{a}\to\mathcal F_{\bar m}^{a}$
is surjective.
Taking $a=0$ shows that $\rho_N^k$ is surjective.
Consequently, $\rho_T^k$ is surjective, and the assertion follows.
\end{proof}

In our case (Condition~\ref{cond:01-condition}), Condition~\ref{cond:odd-intersection-cohomology-vanish}
holds for both $\bar t$ and $\bar m$ by
Theorem~\ref{thm:intersection-cohomology-group},
and the nonequivariant change of perversity $\rho$
is surjective by Proposition~\ref{prop:change-of-perversity}.

\subsection{The quotient ring structure}\label{subsec:equivariant-IH-quotient-ring}
In this subsection, we define the quotient ring structure on $I_{\bar m}H_T^\ast(X)$.
For each $N$, we use the canonical identification
\[
I_{\bar t}H^\ast(X_{T,N})
\cong
H^\ast(X_{T,N})
\]
and equip $I_{\bar t}H^\ast(X_{T,N})$ with the ordinary cup-product
ring structure.

We first show that $\ker\rho_T$ is a homogeneous ideal.

\begin{lemma}
\label{lem:kernel-rho-T-ideal}
Suppose that Condition~\ref{cond:odd-intersection-cohomology-vanish}
holds for both $\bar t$ and $\bar m$.
The kernel of
\[
\rho_T:
I_{\bar t}H_T^\ast(X)
\to
I_{\bar m}H_T^\ast(X)
\]
is a homogeneous ideal of $I_{\bar t}H_T^\ast(X)$.
\end{lemma}

\begin{proof}
Since $\rho_T$ preserves degree, its kernel is a homogeneous subgroup.
It remains to show that it is an ideal.
Let
$a\in I_{\bar t}H_T^i(X)$ and
$x\in\ker\rho_T\cap I_{\bar t}H_T^j(X)$.
Choose $N$ sufficiently large so that the degree-wise stabilization
holds in degrees $i$, $j$, and $i+j$.
Then $a$ and $x$ correspond to classes
\[
a_N\in I_{\bar t}H^i(X_{T,N}),
\qquad
x_N\in I_{\bar t}H^j(X_{T,N}),
\]
and the product $ax$ corresponds to $a_Nx_N$.
Since $x\in\ker\rho_T$, the compatibility of change of perversity with
the stabilization isomorphisms gives
$\rho_N(x_N)=0$.
Thus
$x_N\in\ker\rho_N$.
By the naturality of the intersection cup product
\cite[Proposition~7.3.5]{Fri20}
and its compatibility with the ordinary cup product
\cite[Proposition~7.3.9]{Fri20}, we have
\[
\rho_N(a_Nx_N)
=\rho_N(x_N)\smile_{\bar m,\bar t}a_N=0.
\]
Passing back through the stabilization isomorphism of
Lemma~\ref{lem:degree-wise-stabilization}, we obtain
$\rho_T(ax)=0$.
Therefore $ax\in\ker\rho_T$.
Hence $\ker\rho_T$ is a homogeneous ideal of
$I_{\bar t}H_T^\ast(X)$.
\end{proof}

We can now define the ring structure on middle-perversity equivariant
intersection cohomology.

\begin{definition}
\label{def:ring-structure-of-equivariant-intersection-cohomology}
Under the assumptions of
Lemma~\ref{lem:surjectivity-of-equivariant-change-of-perversity},
the equivariant change of perversity
\[
\rho_T:
H_T^\ast(X)\cong I_{\bar t}H_T^\ast(X)
\to
I_{\bar m}H_T^\ast(X)
\]
is surjective.
Together with Lemma~\ref{lem:kernel-rho-T-ideal}, it induces an
isomorphism of graded abelian groups
\[
H_T^\ast(X)/\ker\rho_T
\xrightarrow{\cong}
I_{\bar m}H_T^\ast(X).
\]
We equip $I_{\bar m}H_T^\ast(X)$ with the ring structure transported
through this isomorphism from the quotient ring
$H_T^\ast(X)/\ker\rho_T$.
We call the resulting graded $H^\ast(BT)$-algebra the
\emph{middle-perversity equivariant intersection cohomology ring} of
$X$.
\end{definition}

\section{Combinatorial description of equivariant intersection cohomology}
\label{sec:combinatorial-description-equivariant-intersection-cohomology}
In this section, we give a combinatorial description of the middle-perversity equivariant intersection cohomology ring. We first construct the contracted labelled graph and then determine the kernel of the equivariant change of perversity.

\subsection{The contracted graph and its face ring}
\label{subsec:contracted-graph-face-ring}

Let $w_1,\ldots,w_r$ be the regular vertices of $P$, listed in cyclic order.
Indices in this section are understood modulo $r$.

For each $j$, let $C_j$ be the path in $\Gamma$ from $w_{j-1}$ to $w_j$ whose internal vertices, if any, are $0$-vertices.
Thus, the paths $C_1,\ldots,C_r$ decompose the cycle $\Gamma$ at its regular vertices.

At every $0$-vertex, the two adjacent characteristic vectors agree up to sign.
Equivalently, the corresponding primitive characters agree up to sign.
Consequently, at every internal vertex of $C_j$, the two adjacent edge labels agree up to sign.
Hence all edge labels on $C_j$ generate the same principal ideal in $H^\ast(BT)$.
Choose a primitive element
\[
\beta_j\in H^2(BT)
\]
generating this common ideal.

\begin{figure}[htbp]
\centering
\begin{tikzpicture}
  \draw (0,0) -- (0.9,0.8)
    node[midway,above left] {$\beta_j$};
  \draw (0.9,0.8) -- (1.9,0.8)
    node[midway,above] {$\pm\beta_j$};
  \draw (1.9,0.8) -- (2.2,0.8);
  \node at (2.6,0.8) {$\cdots$};
  \draw (3.0,0.8) -- (3.3,0.8);
  \draw (3.3,0.8) -- (4.3,0.8)
    node[midway,above] {$\pm\beta_j$};
  \draw (4.3,0.8) -- (5.2,0)
    node[midway,above right] {$\beta_j$};
  \node[fill=white,inner sep=0pt] at (0.9,0.8) {$\circ$};
  \node[fill=white,inner sep=0pt] at (1.9,0.8) {$\circ$};
  \node[fill=white,inner sep=0pt] at (3.3,0.8) {$\circ$};
  \node[fill=white,inner sep=0pt] at (4.3,0.8) {$\circ$};
  \node at (0,0) {$\bullet$};
  \node[below] at (0,-0.1) {$w_{j-1}$};
  \node at (5.2,0) {$\bullet$};
  \node[below] at (5.2,-0.1) {$w_j$};
  \node at (2.6,1.55) {$C_j\subset\Gamma$};
  \draw[->] (5.9,0.4) -- (7.0,0.4);
  \draw (7.7,0) -- (10.5,0)
    node[midway,above] {$\beta_j$};
  \node at (7.7,0) {$\bullet$};
  \node[below] at (7.7,-0.1) {$w_{j-1}$};
  \node at (10.5,0) {$\bullet$};
  \node[below] at (10.5,-0.1) {$w_j$};
  \node at (9.1,1.5)
    {$\widetilde e_j\subset\widetilde\Gamma$};
\end{tikzpicture}
\caption{Replacing $C_j$ by $\widetilde e_j$.
Open circles denote $0$-vertices.
Edge labels are shown up to sign.}
\label{fig:contracted-graph}
\end{figure}
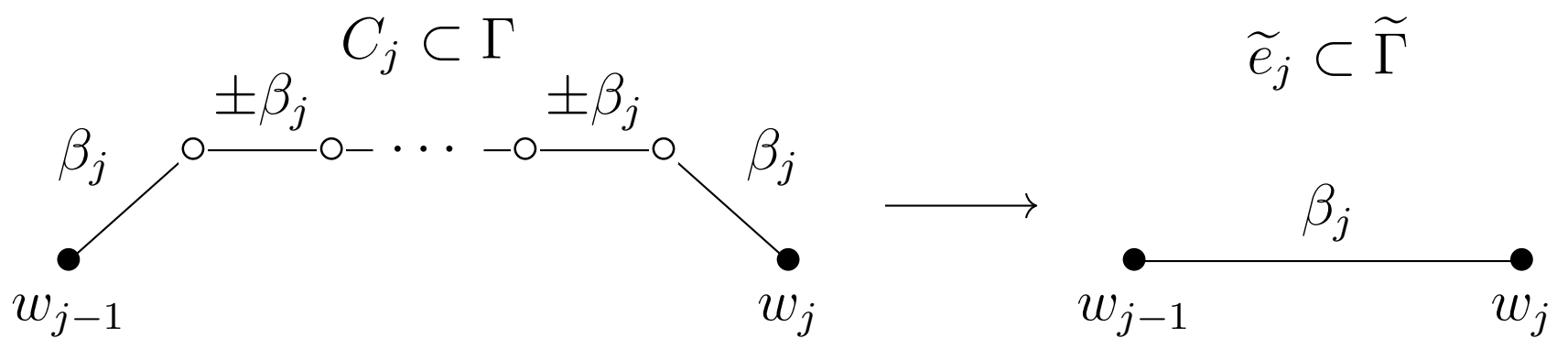

We now define the contracted graph.

\begin{definition}[contracted graph]\label{def:contracted graph}
The \emph{contracted graph} $\widetilde{\Gamma}$ is the labelled graph with vertex set
\[
V(\widetilde{\Gamma})=\{w_1,\ldots,w_r\},
\]
and, for each $j$, an edge $\widetilde e_j$ with endpoints $w_{j-1}$ and $w_j$, labelled by $\beta_j$.
Equivalently, $\widetilde{\Gamma}$ is obtained from $\Gamma$ by contracting every maximal path whose internal vertices are $0$-vertices.
\end{definition}

\begin{remark}
The choice of the sign of each $\beta_j$ is immaterial for the congruence conditions below, since only the principal ideal $(\beta_j)$ is involved.
\end{remark}

The construction of $\widetilde{\Gamma}$ is the graph-theoretic
counterpart of the resolution in
Proposition~\ref{prop:resolution-by-quasitoric}.
Indeed, removing the $0$-vertices and replacing each path $C_j$ by a
single edge gives the labelled graph associated with the nonsingular
part $M$ appearing there.
We use the notion of a torus graph in the sense of \cite{MMP07}.

We first show that the contracted graph is a torus graph.

\begin{lemma}\label{lem:contracted-graph-torus-graph}
The labelled graph $\widetilde{\Gamma}$ admits the structure of a
$2$-valent torus graph.
\end{lemma}

\begin{proof}
Since each $w_j$ is regular, $\beta_j$ and $\beta_{j+1}$
form a basis of $H^2(BT)\cong\mathbb Z^2$.
Hence $\widetilde{\Gamma}$ admits the structure of a
$2$-valent torus graph.
\end{proof}

We next consider the graph equivariant cohomology of the contracted
graph $\widetilde{\Gamma}$, using the same convention as in
Definition~\ref{def:graph-eq-cohomology-of-Gamma}.

\begin{definition}[graph equivariant cohomology]
\label{def:graph-equivariant-cohomology}
The \emph{graph equivariant cohomology} of $\widetilde{\Gamma}$ is defined by
\[
H^{\ast}_{T}(\widetilde{\Gamma})
:=
\left\{
(h_1,\ldots,h_r)\in H^{\ast}(BT)^r
\ \middle|
h_j-h_{j-1}\in(\beta_j)
\text{ for every }j
\right\}.
\]
It is a graded $H^{\ast}(BT)$-algebra under componentwise multiplication.
\end{definition}

By \cite[Theorem~5.5]{MMP07}, the equivariant cohomology of a torus graph is isomorphic to the face ring of its associated simplicial poset. Hence
\[
H^{\ast}_{T}(\widetilde{\Gamma})
\cong
\mathbb Z[\mathcal P(\widetilde{\Gamma})],
\]
where $\mathcal P(\widetilde{\Gamma})$ denotes the simplicial poset associated with the torus graph $\widetilde{\Gamma}$.

\begin{remark}
When $r=2$, the graph $\widetilde{\Gamma}$ consists of two vertices joined by two edges.
In this case, $\mathcal P(\widetilde{\Gamma})$ is a simplicial poset which is not the face poset of a simplicial complex.
When $r\geq3$, $\mathcal P(\widetilde{\Gamma})$ is the face poset of a simplicial complex, and its face ring is the usual Stanley-Reisner ring.
\end{remark}

Recall that $X$ is a $4$-dimensional locally standard
$T^2$-pseudomanifold over a convex polygon $P$ satisfying
Condition~\ref{cond:01-condition}.
By Theorem~\ref{thm:str-of-eq-cohomology},
$H_T^\ast(X)\cong H_T^\ast(\Gamma)$.
Under this identification, we regard an element of $H_T^\ast(X)$
as a function $f:V(\Gamma)\to H^\ast(BT)$.
We now define the restriction map to the regular vertices.

\begin{definition}[restriction map]\label{def:restriction-map}
Define $c:H_T^\ast(X)\to H^\ast(BT)^r$
by
\[
c(f):=
\bigl(f(w_1),\ldots,f(w_r)\bigr).
\]
\end{definition}

Let $z_1,\ldots,z_s$ be the $0$-vertices of $P$.
For each $z_i$, Proposition~\ref{prop:zero-vertex-class} gives a class
\[
\zeta_i\in H_T^2(X)
\]
supported at $z_i$.
More precisely, if $z_i$ is an internal vertex of $C_j$, we choose
the sign of $\zeta_i$ so that
\[
\zeta_i(z_i)=\beta_j,
\qquad
\zeta_i(v)=0
\quad
\text{for every }v\neq z_i.
\]

These classes describe the kernel of $c$
and give the following quotient description.

\begin{proposition}\label{prop:restriction-map}
The restriction map
$c:H_T^\ast(X)\to H_T^\ast(\widetilde{\Gamma})$
is a surjective homomorphism of graded $H^\ast(BT)$-algebras.
Its kernel is
\[
\ker c
=
\bigoplus_{i=1}^{s}H^\ast(BT)\zeta_i
=
(\zeta_1,\ldots,\zeta_s)
\subset H_T^\ast(X).
\]
Consequently, $c$ induces an isomorphism of graded
$H^\ast(BT)$-algebras
\[
H_T^\ast(X)/(\zeta_1,\ldots,\zeta_s)
\cong
H_T^\ast(\widetilde{\Gamma})
\cong
\mathbb Z[\mathcal P(\widetilde{\Gamma})].
\]
\end{proposition}

\begin{proof}
Let $f\in H_T^\ast(X)$.
Since every edge label along $C_j$ generates $(\beta_j)$,
the edge congruences imply
$f(w_j)-f(w_{j-1})\in(\beta_j)$ for every $j$.
Hence $c(f)\in H_T^\ast(\widetilde{\Gamma})$.
Since $c$ is given by restriction to the regular vertices,
it preserves multiplication and the $H^\ast(BT)$-action.
Let $h=(h_1,\ldots,h_r)\in H_T^\ast(\widetilde{\Gamma})$.
Set $f(w_j):=h_j$ for every regular vertex $w_j$.
For each $j$, set $f(v):=h_{j-1}$
at every internal vertex $v$ of $C_j$.
Since $h_j-h_{j-1}\in(\beta_j)$,
the edge congruence relations hold along every $C_j$.
By Theorem~\ref{thm:str-of-eq-cohomology}, this function
determines an element $f\in H_T^\ast(X)$.
By construction, $c(f)=h$.
Therefore $c$ is surjective.

Since $c$ is an algebra homomorphism and $c(\zeta_i)=0$,
we have $(\zeta_1,\ldots,\zeta_s)\subset\ker c$.
Conversely, let $f\in\ker c$.
For each path $C_j$, the edge congruences and
$f(w_{j-1})=0$ imply that $f(z_i)\in(\beta_j)$
at every internal vertex $z_i$.
Thus we can write $f(z_i)=q_i\beta_j$
with $q_i\in H^\ast(BT)$.
By the defining property of $\zeta_i$, we obtain
\[
f=\sum_{i=1}^{s}q_i\zeta_i.
\]
Therefore
$\ker c=\sum_{i=1}^{s}H^\ast(BT)\zeta_i
=(\zeta_1,\ldots,\zeta_s)$.
Since the classes $\zeta_i$ have pairwise disjoint supports,
this sum is direct.
We finally verify the stated isomorphisms.
The first isomorphism follows from the surjectivity of $c$
and the description of its kernel.
The second isomorphism follows from
\cite[Theorem~5.5]{MMP07}.
\end{proof}

\subsection{The kernel of the equivariant change of perversity}
\label{subsec:kernel-equivariant-change-of-perversity}

We first compare the kernel of the equivariant change of perversity
with the kernel of the restriction map $c$.

\begin{lemma}\label{lem:subset-of-ker}
One has $\ker\rho_T\subset\ker c$.
\end{lemma}

\begin{proof}
Let $a\in\ker\rho_T$, and let $w_j$ be a regular vertex.
For each $N$, let
\[
i_{j,N}:B_NT=(w_j)_{T,N}\to X_{T,N}
\]
be the map induced by the inclusion of the fixed point corresponding to
$w_j$.
Here the pullback on intersection cohomology is defined
via stratification independence near the regular fixed point.
Let $a_N\in I_{\bar t}H^\ast(X_{T,N})$ denote the component of $a$ at
the $N$-th finite-dimensional Borel approximation.
By the naturality
of change of perversity, the following diagram commutes:
\[
\begin{tikzcd}[ampersand replacement=\&]
I_{\bar t}H^\ast(X_{T,N})
\ar[r,"\rho_N"]
\ar[d,"i_{j,N}^\ast"']
\&
I_{\bar m}H^\ast(X_{T,N})
\ar[d,"i_{j,N}^\ast"]
\\
I_{\bar t}H^\ast(B_NT)
\ar[r,"\cong"']
\&
I_{\bar m}H^\ast(B_NT).
\end{tikzcd}
\]
Here the bottom horizontal map is an isomorphism because $B_NT$ is a
manifold.
Since $a\in\ker\rho_T$, we have $\rho_N(a_N)=0$ for every $N$.
The commutative diagram therefore gives $i_{j,N}^\ast(a_N)=0$ for every
$N$. Passing to the inverse limit, we obtain
$a|_{w_j}=0$ in $H_T^\ast(w_j)\cong H^\ast(BT)$.
Since this holds for every regular vertex $w_j$, we have $c(a)=0$.
Therefore $a\in\ker c$, and hence $\ker\rho_T\subset\ker c$.
\end{proof}

We now show that the inclusion above is an equality.

\begin{theorem}\label{thm:ker=zeta}
The kernel of the equivariant change of perversity is
\[
\ker\rho_T
=
(\zeta_1,\ldots,\zeta_s).
\]
\end{theorem}

\begin{proof}
By Lemma~\ref{lem:subset-of-ker}
and Proposition~\ref{prop:restriction-map}, we have
\[
\ker\rho_T
\subset
\ker c
=
\bigoplus_{i=1}^s H^\ast(BT)\zeta_i.
\]
It remains to prove the reverse inclusion.
Fix a degree $k$. 
By Lemma~\ref{lem:eq-IH-group-structure},
\begin{align*}
  \rank H_T^k(X)
&=
\sum_{a+b=k}
\rank H^a(BT)\,\rank H^b(X),
\\
\rank I_{\bar m}H_T^k(X)
&=
\sum_{a+b=k}
\rank H^a(BT)\,\rank I_{\bar m}H^b(X).
\end{align*}
By Proposition~\ref{prop:change-of-perversity}, the nonequivariant
change of perversity is surjective and its kernel is
$\bigoplus_{i=1}^s\mathbb Z u_i$, with each $u_i$ of degree $2$.
Hence
\[
\rank H_T^k(X)-\rank I_{\bar m}H_T^k(X)
=
s\,\rank H^{k-2}(BT).
\]
Since $\rho_T^k$ is surjective by
Lemma~\ref{lem:surjectivity-of-equivariant-change-of-perversity}, it follows that
$\rank\ker\rho_T^k=s\,\rank H^{k-2}(BT)$.
On the other hand, Proposition~\ref{prop:restriction-map} gives
\[
(\ker c)^k
=
\bigoplus_{i=1}^s H^{k-2}(BT)\zeta_i,
\]
so $\ker\rho_T^k$ and $(\ker c)^k$ have the same rank. Since
$\ker\rho_T^k\subset(\ker c)^k$, restriction of $\rho_T^k$ gives a short
exact sequence
\[
0
\to
\ker\rho_T^k
\to
(\ker c)^k
\to
\rho_T^k\bigl((\ker c)^k\bigr)
\to
0.
\]
Therefore $\rho_T^k((\ker c)^k)$ has rank zero. 
By Lemma~\ref{lem:eq-IH-group-structure} and
Theorem~\ref{thm:intersection-cohomology-group},
$I_{\bar m}H_T^k(X)$ is free abelian, and hence so is its subgroup
$\rho_T^k((\ker c)^k)$.
It follows that this subgroup is zero.
Thus $(\ker c)^k\subset\ker\rho_T^k$, and hence
$\ker\rho_T^k=(\ker c)^k$.
Since this holds for every $k$,
\[
\ker\rho_T
=
\ker c
=
(\zeta_1,\ldots,\zeta_s).
\]
\end{proof}

We now combine the quotient description with the kernel computation.
This gives the following combinatorial description.

\begin{theorem}[face-ring description]
\label{thm:Stanley-Reisner-description}
With the ring structure of
Definition~\ref{def:ring-structure-of-equivariant-intersection-cohomology},
there are isomorphisms of graded $H^\ast(BT)$-algebras
\[
I_{\bar m}H_T^\ast(X)
\cong
H_T^\ast(X)/\ker\rho_T
=
H_T^\ast(X)/(\zeta_1,\ldots,\zeta_s)
\cong
H_T^\ast(\widetilde{\Gamma})
\cong
\mathbb Z[\mathcal P(\widetilde{\Gamma})].
\]
\end{theorem}

\begin{proof}
The first isomorphism follows from
Definition~\ref{def:ring-structure-of-equivariant-intersection-cohomology},
and the equality follows from Theorem~\ref{thm:ker=zeta}.
The remaining isomorphisms follow from
Proposition~\ref{prop:restriction-map}.
\end{proof}

As a consequence, the middle-perversity equivariant intersection
cohomology of $X$ can be identified with the equivariant cohomology of
the nonsingular part obtained by removing the $0$-vertices.

\begin{corollary}
\label{cor:eq-intersection-cohomology-of-M}
Let $M$ be the nonsingular part of $X$ in
Proposition~\ref{prop:resolution-by-quasitoric}.
Then there is an isomorphism of graded $H^\ast(BT)$-algebras
\[
I_{\bar m}H_T^\ast(X)
\cong
H_T^\ast(M).
\]
\end{corollary}

\begin{proof}
By construction, the labelled graph of $M$ is
$\widetilde{\Gamma}$.
Since $M$ is a torus manifold with vanishing odd-degree cohomology,
\cite{MP06} identifies the image of the fixed-point
restriction map with the graph equivariant cohomology of
$\widetilde{\Gamma}$.
Hence
\[
H_T^\ast(M)\cong H_T^\ast(\widetilde{\Gamma})
\]
as graded $H^\ast(BT)$-algebras.
The assertion follows from
Theorem~\ref{thm:Stanley-Reisner-description}.
\end{proof}

For the basic Thom block, the preceding corollary takes the following form.

\begin{example}
The Thom block $\Th$ has an orbit triangle
with exactly two regular vertices.
Thus $r=2$, so $M=S^4$ in
Corollary~\ref{cor:eq-intersection-cohomology-of-M}.
Hence
\[
I_{\bar m}H_T^\ast(\Th)
\cong
H_T^\ast(S^4)
\]
as graded $H^\ast(BT)$-algebras.
\end{example}

\appendix

\section{Cohomology of connected sums}
\label{sec:appendix-cohomology}

In this appendix, we collect connected-sum formulas for ordinary and intersection cohomology used in the computations of Part~\ref{part-2}.
Throughout this section, $X\#Y$ denotes the ordinary connected sum of
oriented stratified pseudomanifolds formed at regular points.
Here, a regular point means a point of the regular stratum.
We use the connected-sum construction of
\cite[Definition~10.6]{BJ82}.

\subsection{Ordinary cohomology of connected sums}
\label{subsec:cohomology-connected-sums}

We first describe the cohomology ring of a connected sum in terms of the cohomology rings of its summands.

\begin{theorem}\label{thm:connected-sum-for-cohomology-ring}
  Let $X$ and $Y$ be compact, connected, oriented normal stratified pseudomanifolds of dimension $n\geq 2$,
  and let $X \# Y$ be their connected sum at regular points.
Then there is an isomorphism of cohomology rings
  \[
  H^{\ast}(X \# Y) \cong
  H^{\ast}(X \vee Y)/\langle \omega_X - \omega_Y \rangle,
  \]
  where $\omega_X \in H^n(X)$ and $\omega_Y \in H^n(Y)$ are the fundamental cohomology classes.
\end{theorem}

\begin{proof}
The proof follows the same argument as in the manifold case (see
\cite[Lemma~2.1(2)]{BCFHKLN21}).
The argument uses only the local manifold structure near the
connected-sum sphere and the fundamental cohomology classes.
Since the connected sum is formed at regular points, the same argument
applies in the present setting.
Collapsing the separating sphere gives a cofibration sequence
\[
S^{n-1}\to X\#Y\to X\vee Y\to S^n.
\]
As in the manifold case, the induced map
$H^n(S^n)\to H^n(X\vee Y)$ is injective, with image generated by
$\omega_X-\omega_Y$.
The result follows from the associated long exact cohomology sequence.
\end{proof}

\subsection{Intersection cohomology of connected sums}
\label{subsec:intersection-cohomology-connected-sums}

We next give the corresponding additive formula for intersection cohomology.

\begin{theorem}
\label{thm:connected-sum-for-intersection-cohomology}
Let $X$ and $Y$ be compact, connected, oriented normal stratified pseudomanifolds of dimension $n\geq2$, and let $X\#Y$ be their connected sum at regular points.
Then, for any GM perversity $\bar p$,
\[
I_{\bar p}H^i(X\#Y)
\cong
\begin{cases}
\mathbb Z, & (i=0,n),\\
I_{\bar p}H^i(X)\oplus I_{\bar p}H^i(Y), & (0<i<n),\\
0, & (\text{otherwise}).
\end{cases}
\]
\end{theorem}

We first record a lemma describing the effect of removing a ball from the regular stratum.

\begin{lemma}\label{lem:removing-regular-ball}
Let $X$ be a compact, connected, oriented normal stratified pseudomanifold
of dimension $n$, and let $D\subset X$ be a closed $n$-ball contained in
the regular stratum.
Set $X':=X\setminus\mathring D$.
Then, for $i<n$, the inclusion $X'\to X$ induces an isomorphism
\[
I_{\bar p}H^i(X)\xrightarrow{\cong}I_{\bar p}H^i(X'),
\]
and
$I_{\bar p}H^n(X')=0$.
\end{lemma}

\begin{proof}
By excision,
\[
I_{\bar p}H^i(X,X')
\cong H^i(D^n,S^{n-1})
\cong
\begin{cases}
\Z, & (i=n),\\
0, & (i\neq n).
\end{cases}
\]
The map
$H^n(D^n,S^{n-1})\to I_{\bar p}H^n(X)$
sends the relative fundamental cohomology class
to the fundamental cohomology class of $X$.
Since both groups are isomorphic to $\mathbb Z$,
it is an isomorphism.
The assertion follows from the long exact sequence
of the pair $(X,X')$.
\end{proof}

\begin{proof}[Proof of Theorem~\ref{thm:connected-sum-for-intersection-cohomology}]
Since $X\#Y$ is a compact, connected, oriented normal stratified
pseudomanifold of dimension $n$, we have
\[
I_{\bar p}H^0(X\#Y)\cong I_{\bar p}H^n(X\#Y)
\cong\mathbb Z.
\]

Let $D_X\subset X$ and $D_Y\subset Y$ be the closed
$n$-balls used in the connected-sum construction.
Slightly enlarging the images of
$X\setminus\mathring D_X$ and $Y\setminus\mathring D_Y$
across the gluing sphere gives closed neighborhoods $U$ and $V$
whose interiors cover $X\#Y$.
These neighborhoods satisfy
$U\cong X\setminus\mathring D_X$ and
$V\cong Y\setminus\mathring D_Y$.
We may choose them so that
$U\cap V\cong S^{n-1}\times[0,1]$,
with $U\cap V$ contained in the regular stratum.
By the preceding lemma, for $i<n$,
$I_{\bar p}H^i(U)\cong I_{\bar p}H^i(X)$,
$I_{\bar p}H^i(V)\cong I_{\bar p}H^i(Y)$,
and
$I_{\bar p}H^n(U)=I_{\bar p}H^n(V)=0$.
Since $U\cap V$ is contained in the regular stratum,
its intersection cohomology agrees with its ordinary cohomology.
Hence
\[
I_{\bar p}H^i(U\cap V)
=
H^i(U\cap V)
\cong H^i(S^{n-1}).
\]
For $0<i<n$, the Mayer--Vietoris sequence therefore gives
$I_{\bar p}H^i(X\#Y)
\cong I_{\bar p}H^i(X)\oplus I_{\bar p}H^i(Y)$,
as required.
\end{proof}

\begin{remark}\label{rem:connected-sum-change-of-perversity}
The isomorphisms in Theorem~\ref{thm:connected-sum-for-intersection-cohomology} are obtained from excision, long exact sequences of pairs, and Mayer--Vietoris sequences.
Hence, by Remark~\ref{rem:change-of-perversity}, they are compatible with change of perversity.
In particular, for $\bar p\leq\bar q$ and $0<i<n$, the following diagram commutes.
\[
\begin{tikzcd}[ampersand replacement=\&]
I_{\bar q}H^i(X\#Y)
\ar[r,"\cong"]
\ar[d]
\&
I_{\bar q}H^i(X)\oplus I_{\bar q}H^i(Y)
\ar[d]
\\
I_{\bar p}H^i(X\#Y)
\ar[r,"\cong"']
\&
I_{\bar p}H^i(X)\oplus I_{\bar p}H^i(Y)
\end{tikzcd}
\]
\end{remark}

\section*{Acknowledgments}
The authors would like to thank Jongbaek Song
for telling us about his paper.

\section*{Declaration of generative AI and AI-assisted technologies in the manuscript preparation process}
During the preparation of this work, the authors used ChatGPT (OpenAI, GPT-5.6 Sol) in order to improve the language, clarity, and organization of the manuscript. After using this tool, the authors reviewed and edited the content as needed and take full responsibility for the content of the published article.

\end{document}